\documentclass[aos]{imsart}

\RequirePackage{amsthm,amsmath,amsfonts,amssymb}
\RequirePackage[numbers]{natbib}
\RequirePackage[colorlinks,citecolor=blue,urlcolor=blue]{hyperref}
\RequirePackage{graphicx}

\startlocaldefs
\usepackage[ruled,linesnumbered]{algorithm2e}

\usepackage{xr}
\usepackage{arydshln}
\usepackage{ulem}
\usepackage{mathrsfs}
\usepackage{algpseudocode}
\usepackage{enumerate}
\usepackage{bm}

\newtheorem{lemma}{Lemma}

\def\SM{\textit{SM}}

\usepackage{tikz}
\usetikzlibrary{shapes,arrows}

\tikzstyle{process} = [rectangle,draw=black!70,fill=white!30,rounded corners]
\tikzstyle{arrow} = [color=black,thick]

\numberwithin{equation}{section}

\newtheorem{theorem}{Theorem} 

\newtheorem{remark}{Remark} 
\numberwithin{equation}{section}

\DeclareMathOperator*{\argmin}{arg\,min}
\DeclareMathOperator*{\argmax}{arg\,max}

\usepackage{multirow}
\def\E{{\mathbf E}}
\def\P{{\mathbf P}}

\def\D{{\mathcal D}}

\def\B{{\mathbb{B}}}
\def\R{{\mathbb{R}}}
\def\A{{\small\mathcal{A}}}
\def\C{{\small\mathcal{C}}}

\allowdisplaybreaks

\endlocaldefs

\begin{document}

\begin{frontmatter}
\title{Ensemble Projection Pursuit  for General Nonparametric Regression}
\runtitle{Ensemble Projection Pursuit}

\begin{aug}
\author{\fnms{Haoran}~\snm{Zhan}\ead[label=e1]{haoran.zhan@u.nus.edu}},
\author{\fnms{Mingke}~\snm{Zhang}\ead[label=e2]{mingke.zhang@u.nus.edu}}
\and
\author{\fnms{Yingcun}~\snm{Xia}\ead[label=e3]{yingcun.xia@nus.edu.sg}}
\address{\normalsize Department of Statistics and Data Science, National University of Singapore}

\address{\printead[presep={ \ }]{e1} \printead[presep={, \ }]{e2} \printead[presep={, \ }]{e3}}

\end{aug}

\begin{abstract}


The projection pursuit regression (PPR) has played an important role in the development of statistics and machine learning. However, when compared to other established methods like random forests (RF) and support vector machines (SVM), PPR has yet to showcase a similar level of accuracy as a statistical learning technique. In this paper, we revisit the estimation of PPR and propose an \textit{optimal} greedy algorithm and an ensemble approach via "feature bagging", hereafter referred to as ePPR, aiming to improve the efficacy. Compared to RF, ePPR has two main advantages. Firstly, its theoretical consistency can be proved for more general regression functions as long as they are $L^2$ integrable, and higher consistency rates can be achieved. Secondly,  ePPR does not split the samples, and thus each term of PPR is estimated using the whole data, making the minimization more efficient and guaranteeing the smoothness of the estimator.  Extensive comparisons based on real data sets show that ePPR is  more efficient in regression and classification than RF and other competitors. The efficacy of ePPR, a variant of Artificial Neural Networks (ANN), demonstrates that with suitable statistical tuning, ANN can equal or even exceed RF in dealing with small to medium-sized datasets. This revelation challenges the widespread belief that ANN's superiority over RF is limited to processing extensive sample sizes.


\end{abstract}



\end{frontmatter}

\section{Introduction}



The technique of projection pursuit, which was first introduced by \cite{kruskal1969toward} and \cite{Friedman1974}, aims to reveal the most significant structures within high-dimensional data by sequentially
utilizing lower-dimensional projections.  Once a projection has been identified, the structure that lies along it is removed, and the procedure is then reiterated to detect the next projection.
This idea was later implemented for regression  by \cite{friedman1981projection}, who referred to this technique as the projection pursuit regression (PPR) and popularized it in the 1980s.
 PPR is well-suited for estimating general nonparametric regressions due to its capability to approximate any general function \citep{Lee1992, petrushev1998approximation, Rashid2020}.

More specifically, consider the general regression problem with response $ Y $ and predictor $ X \in \R^p$. Our main interest lies in the conditional mean,  $m(X)=\E(Y|X)$.
Following the categorization of two cultures of \cite{Breiman2001}, the projection pursuit regression can be used in two ways. The first way to use PPR is by assuming model $  m(X) =  g_1(\theta_1^\top X) + ... + g_K(\theta_K^\top X)
$, where $K$ is assumed to be a fixed and finite, yet unknown number.  See, for example, \cite{chen1991estimation},  
\cite{HALL1993} and  \cite{xia1999extended}.  It is evident that this statistical model may not be appropriate for real-world data and is not recommended by  \cite{Breiman2001}. The second way to use PPR is by approximating or estimating a general nonparametric regression based on the following fact: for any smooth function $ m $, there exists a sequence of projection directions $ \{\theta_k \in \mathbb{R}^p: k=1, 2,...\} $ and a sequence of ridge functions $\{g_k(.), k = 1, 2, ...\} $ satisfying
\begin{equation}
\mbox{PPR expansion:} \ \ \  m(X) =  \sum_{k=1}^\infty g_k(\theta_k^\top X)\ \ \ a.s..  \label{PPRidentity}
\end{equation}
Under mild conditions, \cite{petrushev1998approximation} and \cite{Rashid2020} established that the PPR expansion holds for a general function $m(x)$ in the $L^2$ sense, and the PPR algorithm proposed in \cite{friedman1981projection} is an implementation of this expansion.

The PPR expansion can be estimated with two different approaches:  \textit{the global minimization} and \textit{the greedy algorithm}. The global minimization is to minimize the loss function over all the ridge functions and projection directions, i.e.,
\begin{equation*}\label{the global minimization supp}
    \min_{g_k, \theta_k, k = 1, 2,... } \E\left( Y - \sum_{k=1}^\infty g_k(\theta_k^\top X)\right)^2.
\end{equation*}
This is the classical estimation method in statistics and machine learning, including the estimation of artificial neural networks (ANN). In practice, the number of terms in the expansion must be finite and goes to infinity with sample size. With i.i.d. data $\{(X_i,Y_i)\}_{i=1}^n$,   the estimator based on the global minimization is  $\hat{m}_n(X)=\sum_{k=1}^{K_n} \hat{g}_k(\hat{\theta}_k^\top X)$, where
\begin{equation*}\label{Theoreticalestimatorabs}
      (\hat{g}_k,\hat{\theta}_k)  = \argmin_{g_k\in V_n,\theta_k\in \Theta^{p}} \frac{1}{n}\sum_{i=1}^n \Big( Y_i-  \sum_{k=1}^{K_n} g_k(\theta_k^\top X_i) \Big)^2
\end{equation*}
and $V_n$ is an approximation function space which will be discussed later. Under some regularity conditions, we can show that $\hat{m}_n$ achieves a nearly optimal consistency rate, i.e.,
$$\E(\hat{m}_n(X)-m(X))^2\leq c(\ln^4/n)^{-{2r}/(2r+p)},$$
where $c>0$ and $r$ is the degree of smoothness of  $m$; see SUPPLEMENTARY MATERIALS (\SM) for the proof.
This method has good statistical asymptotic properties, but in practice, its advantages are hardly observed, perhaps because of the large number of parameters required for estimation.

Instead of estimating $g_k(\theta_k^\top x), k = 1, 2, ...$ together, \textit{the greedy algorithm} estimates  the PPR expansion term by term in a sequential manner, which was originally introduced by \cite{friedman1981projection} in their PPR estimation procedure, as described below.
Let
\begin{equation}
\theta_1 = arg \min_{ \theta: ||\theta||_2=1} \E(Y-\E(Y|X^\top \theta))^2 \label{ippr.1}
\end{equation}
and $ g_1(v)  = \E(Y|X^\top \theta_1 = v) $. For $ k \ge 2 $, let the residual $ \delta_k = Y- \sum_{\iota=1}^{k-1} g_\iota(X^\top \theta_\iota) $, and  calculate
\begin{equation}
\theta_k = arg \min_{ \theta: ||\theta||_2=1} \E( \delta_k -\E(\delta_k |X^\top \theta))^2  \label{ippr.2}
\end{equation}
and $ g_k(v)  = \E(Y|X^\top \theta_k = v) $.  \cite{jones1987conjecture} proved that the PPR algorithm converges in the sense that
\begin{equation}\label{ippr.3}
     \E\Big[ \Big|m(X) - \sum_{k=1}^{K} g_k(X^\top \theta_k)\Big|^2\Big] \to 0, \ \  as \ K \to \infty.
\end{equation}
In the literature, there are various proposed improvements to the aforementioned greedy algorithm. For instance, \cite{Lee1992} proposed the relaxed greedy algorithm (RGA) as a means of accelerating the convergence rate of the PPR algorithm. In RGA, the approximation obtained in step $k$ takes the form of $\alpha\cdot\sum_{\iota=1}^{k-1}{g_\iota(X^T\theta_\iota)}+(1-\alpha)\cdot g_k(X^T\theta_\iota)$ for some $\alpha\in (0,1)$. The orthogonal greedy algorithm (OGA), which is well-known in statistical learning and signal processing (see \cite{cai2011orthogonal} and \cite{barron2008approximation}), is another improvement. In OGA, once $g_k(X^T\theta_k)$ is obtained in the $k$-th iteration, the approximation is updated by projecting $m(X)$ onto the linear space spanned by $\{g_\iota(X^T\theta_\iota)\}_{\iota=1}^k$.



Although PPR has good algorithmic and mathematical guarantees, it has not gained popularity among researchers, nor has it been extensively studied. Existing computer packages of PPR do not provide good prediction accuracy in comparison with the  other popular methods such as RF and the support vector machine  \citep[SVM]{cortes1995support}. \cite{hastie01} suggested that PPR's limited popularity was due to its computational complexity during the 1980s, when computational power was significantly lower than it is today. As a result, the artificial neural network (ANN) became a more popular choice although they share some similarities. Nevertheless, with the  tremendous
developments in computational power and algorithms,
there should be a growing interest in studying PPR more comprehensively. This is the primary motivation behind this paper.


The purpose of this paper is to investigate the greedy algorithm of PPR approximation and its convergence theory. To further improve the prediction accuracy of PPR, we  propose a boosting method, called ensemble PPR (ePPR) hereafter. The fundamental idea is to search for projections in \eqref{ippr.1} and \eqref{ippr.2} in a randomly selected subset of $ X $, in a similar way as that RF selects split variables \citep{breiman2001random}, which is also called ``feature bagging'' \citep{ho1995random}. Thus, each run of the algorithm produces a random PPR approximation of the unknown regression function. The ensemble technique is then employed to generate a final estimate of the regression function.  


\subsection{Related work and our contribution}


Our work includes two parts: the greedy algorithm for PPR and a boosting method for a statistical implementation of PPR. Our contribution is summarized as follows.

\begin{itemize}
    \item An \textit{optimal} greedy algorithm, called the Additive Greedy Algorithm (AGA). The algorithm is  motivated by the orthogonal greedy algorithm \cite[OGA]{pati1993orthogonal} and the additive structure of the expansion.  We will show that AGA achieves the optimal convergence rate of greedy algorithms and is faster than the existing algorithms for PPR expansion (see \cite{barron2008approximation}).

\item An ensemble estimation, called ensemble PPR (ePPR). This is partly motivated by RF of \cite{Breiman2001}; see also \cite{Zhu2015} and \cite{athey2019generalized} for the recent development of RF. ePPR shares one of the core elements of RF: the former uses recursive projection and the latter recursive partitioning. ePPR also  borrows two of the other core elements of RF: random selection of predictors, also known as feature bagging, and averaging the individual estimates for the ensemble.

\item Advantages over RF, the  benchmark of nonparametric regression based on the greedy algorithm. These advantages include: (i)  As noted by \cite{friedberg2020local}, RF has limited ability to fit smooth signals because of the step-functions inherited from the classification and regression tree (CART). In contrast, PPR does not split the samples and uses the whole data in each step, making the minimization more efficient and guaranteeing the estimator smooth. (ii) Moreover, while RF theoretically only works for some special classes of regression functions (see, for example, \cite{scornet2015consistency} and \cite{chi2022asymptotic}),   PPR or ePPR can handle more general regression functions as long as the underlying function is $L^2$ integrable. (iii) Additionally,  ePPR has a higher consistency rate than RF, as demonstrated in
 \cite{klusowski2022large} for additive model  and \cite{chi2022asymptotic} for  slightly more general regression functions.

\item Our work also sheds light on the long debate about when ANN and RF have advantages over each other. As a special variant of ANN, ePPR suggests that  ANN can perform as well as or even better than RF for data of small to medium size when statistical tuning is incorporated. This contradicts the commonly received understanding of ANN and RF that ANN can only have  advantages for data with large sample size.

\end{itemize}



\subsection{Organization of this paper} The rest of this paper is organized as follows. Section \ref{s.2.aga}  introduces some notations, and proposes the  greedy algorithm, AGA, and studies its convergence rate. Section \ref{est. via AGA} considers the consistency of statistical estimators of AGA.
Details of ePPR are presented in  Section \ref{ssestimator}; its statistical consistency  is also studied in the section. Extensive numerical studies of  real data sets, including comparison with the popular machine learning methods such as RF and SVM,  are presented in Section \ref{secReal}.  The technical proofs are provided in Section \ref{SecProof}, while some auxiliary results and lemmas  used in the proofs are given in SUPPLEMENTARY MATERIALS (\SM) in a separate file.
A brief discussion about the problems with the proposed ePPR is given in Section \ref{secConculsion}.

\section{PPR with greedy algorithm}  \label{s.2.aga}

A greedy algorithm always associates with a dictionary $Dic\subseteq \mathscr{H}$, where $\mathscr{H}$ is a Banach space equipped with norm $\|\cdot\|$.  For any dictionary $Dic$,  the following convex hull
$$
B_1(Dic):= \overline{\left\{\sum_{j=1}^m{a_jd_j}: d_j\in Dic, \sum_{j=1}^m{|a_j|}\leq 1,m\in\mathbb{Z}^+\right\}}
$$
is a function class generated by $Dic$,  which is commonly used to analyze the convergence rate of  the greedy algorithm (see \cite{barron2008approximation} and \cite{siegel2022optimal}).
Then, by scaling $B_1(Dic)$, we have a linear sub-space of $\mathscr{H}$,
\begin{equation}\label{scaleddistionary}
       [\mathscr{H}]_1(Dic):=
   \bigcup_{c> 0}{c \cdot B_1(Dic)},
\end{equation}
where $c \cdot B_1(Dic):=\{c\cdot \phi: \  \phi\in B_1(Dic)\}$. Therefore, for any $f\in[\mathscr{H}]_1(Dic)$, the distance between $f$ and $B_1(Dic)$ is equal to
$$
\|f\|_{[\mathscr{H}]_1(Dic)}= \inf\{ c>0: f\in c\cdot B_1(Dic)\}.
$$
 \cite{siegel2021sharp} proved that  $[\mathscr{H}]_1(Dic)$ is also a Banach space w.r.t. norm $\|f\|_{[\mathscr{H}]_1(Dic)}$ if $\sup_{d\in Dic}{\|d\|}<\infty$. Note that $\|f\|_{[\mathscr{H}]_1(Dic)}$  depends on $Dic$ and is called the variation norm or the atomic norm in the literature.

Note that if we use  activation functions to approximate the ridge functions in PPR,  some of  the activation functions share the same direction,  thus the dictionary associated with AGA is
\begin{equation}\label{bdkjASZNdjknad;}
     Dic^s_n:=\left\{  \sum_{j=1}^{J_n}{a_j\cdot\sigma_s(\theta^\top x+b_j): \theta\in\Theta^p, b_j\in[-2,2], \ \sum_{j=1}^{J_n}|a_j|\leq 1}\right\},
\end{equation}
where $\sigma_s(.) $ is an activation function and will be discussed later, while the subscript $ n $ is used to indicate that the choice of $J_n $ depends on the sample size $n$.  Since $Dic^s_n$ is designed to approximate $m(x),x\in\mathbb{B}_1^p$ with bounded support, we will show that a bounded range of those $b_j$ in \eqref{bdkjASZNdjknad;} is sufficient, including [-2, 2] in \eqref{bdkjASZNdjknad;}. Meanwhile, the restriction $\sum_{j=1}^{J_n}|a_j|\leq 1$ is to normalize those possible searching directions in $[\mathscr{H}]_1(Dic_n^s)$. In this case, let $\mathscr{H}=L^2(\mathbb{B}_1^p,X):=\{f: \E_X(f^2(X))<\infty\}$.  Then,  it is known that $L^2(\mathbb{B}_1^p,X)$ is  a Hilbert space equipped with norm $\|f\|:=\left( \E_X(f^2(X))\right)^{1/2}$ and inner product $\langle f,g\rangle:= (\|f+g\|^2-\|f-g\|^2)/4$ for any $f,g\in L^2(\mathbb{B}_1^p,X)$. With the above notations, our AGA in population version can be stated as follows.

\vspace{0.3cm}

{\bf  PPR.aga (population version)}  Starting with iteration $ k = 0 $, let $ \delta_0 = m $ and $ m_0^{aga} = 0 $. For each $k \ge 1$,
calculate
\begin{equation}\label{GBJHBKJNKJqqqr}
\sum_{j=1}^{J_n}{a_{k,j}\cdot\sigma_s(\theta_k^\top X+b_{k,j})}= arg \max_{ d\in Dic_n^s} {|\langle \delta_{k-1}, d\rangle|},
\end{equation}
where $ |\langle \delta_{k-1}, d\rangle| $ is defined above and measures the similarity between $\delta_{k-1}$ and $d\in Dic_n^s$. Another method to find the above parameters is by residual approximation:
\begin{equation}\label{GBJHBKJNKJqqqr''}
    \sum_{j=1}^{J_n}{a_{k,j}\cdot\sigma_s(\theta_k^\top X+b_{k,j})} = arg \min_{ h\in[L^2(\mathbb{B}_1^p,X)](Dic_n^s) } {\left\|\delta_{k-1}-h\right\|}.
\end{equation}
Then, update the approximation by
$$
 m_k^{aga}  = \displaystyle\sum_{j=1}^{J_n}{a'_{1,j}\cdot\sigma_s(\theta_1^\top X+b_{1,j})} + ... + \sum_{j=1}^{J_n}{a'_{k,j}\cdot\sigma_s(\theta_k^\top X+b_{k,j})}
 $$
 with coefficients $ a_{\iota,j}' $s calculated by
\begin{equation}\label{BJKBKJNppppqqqq}
       \argmin_{\substack{a'_{\iota,j}\\ \iota=1,\ldots,k\\ j=1,\ldots, J_n}} \Big\| m(X) - \sum_{j=1}^{J_n}{a'_{1,j}\cdot\sigma_s(\theta_1^\top X+b_{1,j})} - ... - \sum_{j=1}^{J_n}{a'_{k,j}\cdot\sigma_s(\theta_k^\top X+b_{k,j})}\Big\|,
\end{equation}
where directions $(\theta_1, b_1), \ldots, (\theta_k,b_k)$ are fixed in the minimization.  Meanwhile, the residual is updated by $\delta_k=m-m_k^{aga}$.


The AGA approximation, denoted by $m_k^{aga}$, bears resemblance to a single-layer artificial neural network (ANN). However, in AGA each group of $J_n$ neurons shares with the same $\theta$, whereas in ANN, all $\theta$'s are allowed to be different. Based on approximation results for ridge functions and ANN in \cite{petrushev1998approximation} and \cite{yarotsky2017error}, this difference allows for a reduction in the number of $\theta$ values required for the PPR approximation process. Meanwhile, our AGA offers greater flexibility compared to OGA because it optimizes all coefficients ${a_{\iota,j}'}$ in the projection step, as demonstrated in  \eqref{BJKBKJNppppqqqq}.

\begin{theorem}[{upper bound of AGA}]\label{Aga population version}
Let $m_k^{aga}$ be the output of AGA  after $k$-th iteration where the input function is $m\in L^2(\mathbb{B}_1^p,X)$. Suppose the density of $X$ satisfies $c^{-1}<f_X(x)<c$ for some constant $c>0$. Then, we have
\begin{align}
\E_{ X}(m_k^{aga}(X)-m(X))^2\leq &\inf_{h\in [L^2(\mathbb{B}_1^p,X)]_{1}(Dic^s_{n})}\Big\{\E(m(X)-h(X))^2 \nonumber \\
   &+ c(s,p)\left( \|h\|_{[L^2(\mathbb{B}_1^p,X)]_{1}(Dic^s_{n})}\right)^2\cdot k^{-1-\frac{2s+1}{2p}}\Big\},\label{dhjkzdbqqq'}
\end{align}
for each iteration $k\in \mathbb{Z}^+$ and some $c(s,p)>0$, where $[L^2(\mathbb{B}_1^p,X)]_{1}(Dic^s_{n}) $ is defined in the same way as $  [\mathscr{H}]_1(Dic)$.
\end{theorem}

\begin{remark}  \rm
    The first term of the upper bound in  \eqref{dhjkzdbqqq'}  is the distance between $m(x)$ and $B_1(Dic^s_n)$, the function space generated by $Dic_n^s$, while the second term  in  \eqref{dhjkzdbqqq'} indicates that the rate of AGA is $k^{-1-\frac{2s+1}{2p}}$ in the approximation space $B_1(Dic_s^n)$.  We will show  later that the above $k^{-1-\frac{2s+1}{2p}}$ is optimal for global approximation using any linear combination of  $k$ terms in the dictionary $Dic^s_n$.
The constant $c(s,p)$ in  \eqref{dhjkzdbqqq'} is universal. Namely, it does not depend on  $n\in\mathbb{Z}^+$. This fact is essential for the theoretical proof.  If we use  \eqref{GBJHBKJNKJqqqr''} to obtain $(\theta_k, b_{k,1},\ldots,b_{k,J_n}, a_{k,1},\ldots, a_{k,J_n})$, we can show a lower approximation rate $k^{-1}$, which is slightly slower  than the optimal convergence rate.
\end{remark}

We end this section by giving the lower bound of AGA below, from which we can see the rate of AGA in Theorem \ref{Aga population version} is optimal.



\begin{theorem}[{lower bound of AGA}]\label{Lower bound of Aga}
If the density of $X$ satisfies $c^{-1}<f_X(x)<c$ for some constant $c>0$, then for any $Dic_n^s$ and $\alpha'(s,p)>\frac{2s+1}{4p}$ we have
$$
  \sup_{m\in [L^2(\mathbb{B}_1^p,X)]_{1}(Dic_n^s)}\sup_{k\geq 1}{k^{-\frac{1}{2}-\alpha'(s,p)}\|m(X)-m^{aga}_k(X)\|_{L^2(\mathbb{B}_1^p, X)}} = \infty.
$$
\end{theorem}

From  \eqref{optimal aga rate population bb} in Theorem \ref{Aga population version},  we have the following convergence result:
\begin{equation}\label{HJKSBKJBKNKJNKNqs}
  \sup_{m\in [L^2(\mathbb{B}_1^p,X)]_{1}(Dic_n^s)}\sup_{k\geq 1}{k^{-\frac{1}{2}-\frac{2s+1}{4p}}\|m(X)-m^{aga}_k(X)\|_{L^2(\mathbb{B}_1^p, X)}} < \infty,
\end{equation}
which indicates that AGA  has an approximation rate at least $\gamma(s,p):=1+(2s+1)/{2p}$ for any  function in  class $[L^2(\mathbb{B}_1^p,X)]_{1}(Dic_n^s)$. Theorem \ref{Lower bound of Aga} further implies that $\gamma(s,p)$ is optimal and can not be improved.  In the proof of Theorem \ref{Lower bound of Aga}, we even show that  $\gamma(s,p)$ is the optimal approximation rate  for approximating functions in $[L^2(\mathbb{B}_1^p,X)]_{1}(Dic_n^s)$ by $\Sigma_k(Dic_n^s):=\{ \sum_{j=1}^k{\alpha_jd_{n,j}}: d_{n,j}\in Dic^s_{n}, \alpha_j\in\mathbb{R} \}$. In conclusion, our AGA can approximate functions in $[L^2(\mathbb{B}_1^p,X)]_{1}(Dic_n^s)$ with the optimal rate $\gamma(s,p)=1+(2s+1)/{2p}$.

\section{Statistical estimator of PPR expansion via AGA}\label{est. via AGA}

Next, we  apply AGA to the estimation of $m(X)$ given i.i.d. data $\D_n=\{(X_i,Y_i)\}_{i=1}^n$. Denote by  $\mu(x)$  the distribution of $X$. Let $\mathfrak{X}_n$ be the Banach space satisfying
$$
  \|f\|_n:=\left(\frac{1}{n}\sum_{i=1}^n{f^2(X_i)}\right)^{\frac{1}{2}}<\infty
$$
for any function $f\in \mathfrak{X}_n$, which is defined  on sample points $(X_1,\ldots,X_n)$. Then, we know $\mathfrak{X}_n$  is also a Hilbert space equipped with an inner product
$\langle f,g\rangle_n=\frac{1}{n}\sum_{i=1}^n{\left(f(X_i)\cdot g(X_i)\right)}$ for any real functions $f,g\in \mathfrak{X}_n$. Note that $\mathbb{Y}=(Y_1,\ldots,Y_n)^\top\in \mathfrak{X}_n$, since it can be regarded as a function on $(X_1,\ldots,X_n)$. With these notations, our sample version of AGA is as follows.

\vspace{0.3cm}

{\bf  PPR.aga (sample version)}
In the sample version of AGA, the input function is changed to $\mathbb{Y}$, defined on $(X_1,\ldots,X_n)$, and $\mathscr{H}=\mathfrak{X}_n$ with empirical norm $\|\cdot\|_n$ and empirical inner product $\langle \cdot,\cdot\rangle_n$.  Then, we only need to replace $\langle \cdot,\cdot\rangle$ in \eqref{GBJHBKJNKJqqqr} and $\|\cdot\|$ in \eqref{BJKBKJNppppqqqq} by their empirical versions, $\langle \cdot,\cdot\rangle_n$ and $\|\cdot\|_n$, respectively.
After $k$ iterations, the output function is denoted by $m^{aga}_k$.


We assume $m(x),x\in\B_1^p$, follows the PPR expansion
\begin{equation}\label{SAjdaldjad..}
    \mathscr{M}_{r,p}:= \left\{ \sum_{\iota=1}^\infty{m_\iota(\theta_\iota^Tx): x\in\B_1^p, \theta_\iota\in\Theta^p, m_\iota\in TV_r[-1,1]}  \right\},
\end{equation}
where $TV_r[-1,1]$ consists of functions defined on $[-1,1]$ whose $r$-th ($r\geq 1$) derivative has bounded total variation (TV).  For any $f(x), x\in [-1,1]$, we denote its total variation and supremum  norm by  $TV(f)$ and $\|f\|_\infty$ respectively. Thus, we also have $TV_r[-1,1]=\{ f: TV(f^{(r)})<\infty\}$. In this case, the consistency rate for truncated AGA estimator $\hat{m}_{k_n}^{aga}:=\max{\{\min{\{m_{k_n}^{aga},t_n\}},-t_n\}}$ is as follows. 

\begin{theorem}[consistency rate]\label{PPRpprmodelrate22}
Assume
$
\E(e^{c\cdot Y^2})<\infty
$
for some $c>0$ and  $m\in\mathscr{M}_{s, p}$  with $\sum_{\iota=1}^\infty {\|m_\iota\|_\infty}<\infty$ and $\sum_{\iota=1}^\infty {TV(m_\iota^{(s)})}<\infty$.  If $t_n\asymp \ln n$, $k_n\asymp n^{(1+\frac{2s-1}{2q})/(2s)}$ and $J_n\asymp n^{(1+\frac{2s-1}{2p})^2/(4s^2)}$,  we have
$$
 \E_{\D_n}\int|\hat{m}^{aga}_{k_n}(x)-m(x)|^2d\mu(x) = O(
n^{-\alpha_*(s,p)}\cdot \ln^4 n),
$$
where  rate $\alpha_*(s,p)= (\frac{s+.5}{p}+1)/\left((\frac{s+.5}{p}+1)+ (1+\frac{1}{2s})(1-\frac{1}{2p})\right)$.
\end{theorem}

\begin{remark} \rm
In this theorem, we only give the statistical consistency rate of the estimator based on AGA when $m\in\mathscr{M}_{s,p}$. In fact,  we can also show its consistency holds for
a more general space than $\mathscr{M}_{s,p}$, for example, one can only assume $m$ is continuous. We do not go into the  details, but a similar proof for the consistency rate of ensemble AGA (Theorem \ref{EPPRresultppr}) will be shown in section \ref{SecProof}. Furthermore, our rate $\alpha_*(s,p)$ is \textit{free from the distribution of} $X$.
\end{remark}

\section{Ensemble PPR estimator by AGA} \label{ssestimator}

To boost the estimation efficiency of AGA, we further use the ``feature bagging'' of \cite{ho1995random} as RF does. Note that at each node of a tree, RF selects the split variable from a random subset of the predictors rather than the full set of predictors $ X $ as CART does \citep{BreiFrieStonOlsh84}. In the same manner, for each projection of \eqref{ippr.1} and \eqref{ippr.2} of PPR, we will randomly select a subset of the predictors for the projection. Because of the random selection, each run of the resultant algorithm  will generate different estimations and predictions. We can then average those  predictions to make a final estimation or prediction. We call this procedure ensemble PPR, denoted by ePPR.

In this section, we start with implementing the PPR algorithm with the projection in a randomly selected subset of $X $, called random PPR (rPPR) hereafter. We can then apply the proposed AGA to rPPR, and call the algorithm random PPR.aga, or rPPR.aga.

For ease of exposition, we introduce several notations that are often used in the rest of this paper. Denote by $\mathcal{A}$ a (random) subset of $ \{ 1, ..., p\} $,  and $\mathcal{A}_\tau,\ \tau=1, 2, ...$ is a sequence of such subsets that may differ from one another. Let $ X_\A $ be a sub-vector of $ X $ consisting of coordinates indexed by elements of $\mathcal{A}$. Similarly, $ x_\A $ can be defined for vector $ x $, and $ X_{i, \A} $ for a sample $ X_i $ of $ X $. To simplify the notation, by writing $ \theta^\top X_\A $ we imply that $ \theta $ has the same dimension as $ X_\A $.  Let $\mathbf{\Omega}_q$ be the collection of all $\A$s with $card(\A) = q$. We can also think of $\mathbf{\Omega}_q $ as a probability space with probability equally assigned to every element. Thus, $ \A_\tau, \tau\geq 1 $, can be regarded as a (random) sample of $ \mathbf{\Omega}_q $. Let $\Xi_k=(\A_1,\ldots, \A_k)$  with  $k\geq 1$, which can be regarded as a random element in the product of probability spaces, $ \mathbf{\Omega}_q ^{\otimes k} $.

\vspace{0.3cm}

{\bf  ePPR.aga} \
We independently execute the sample version of AGA for $B$ times and obtain estimators $m_{k,b}, b = 1, ..., B$. To find each projection and ridge function in an estimator $m_{k,b}$, instead of using all the predictors (or features) in each iteration $ k $,
we randomly select  $\A_{b,k} \subseteq \{1,\ldots,p\}$ and then
 change $Dic_n^s$ to \begin{equation*}
     Dic^{s,q}_n(\A_{b,k}):=\left\{  \sum_{j=1}^{J_n}{a_j\sigma_s(\theta^\top x_{\A_{b,k}}+b_j):  \theta\in\Theta^q, b_j\in[-2,2], \sum_{j=1}^{J_n}|a_j|\leq 1}
     \right\}
\end{equation*} and denote the resulting estimator as  $m_{k,b}$. After running AGA for $B$ times, the ensemble estimator of $m$ is $m_{k}^{ens} = B^{-1} \sum_{b=1}^B m_{k,b}$.  Later, we also write $ k$ as $ k_n$ to indicate it depends on the sample size $ n$.


To make the algorithm more efficient, we can choose   the random subset $ {\cal A}_{b,k}$  that has the maximum similarity in \eqref{GBJHBKJNKJqqqr} or the best approximation in \eqref{GBJHBKJNKJqqqr''}  from $ \ell $ candidates.  Therefore, there are four tuning parameters in ePPR.aga: the number of variables, $ q $, in each random subset for feature bagging, the number of candidate subsets of variables, $ \ell $,   and the number of iterations (or ridges),  $k$, of each individual run of the algorithm, and the number of runs, $ B $, for the ensemble. We have discussed the selection of $ q $ for rPPR.aga, but for the sample version, we also need to consider the estimation efficiency of the projection $ \theta $. Because of the nonparametric estimation of the projections, $ q $ should be relatively much smaller than the sample size $ n $ to get a good estimator of the projection; see \cite{HALL1993}. Thus, we select $ q = \min([2p/3], [n^{0.4}]) $ in our calculation.  For a similar reason, the choice of $ \ell $ also needs to consider the diversity of the estimators $  m_{1,k}, ...,  m_{B,k}$ and the efficiency of each individual $ m_{b,k} $. As an extreme case, if $ \ell = {p \choose q}$ then all the estimators are the same with no diversity. Our limited numerical study suggests that $ \ell \le p/q $ can be used. The selection of $ k $ depends on the model structure, and will be discussed in the next section. For $ B$, there are studies about its selection for RF; see, for example, \cite{Zhang2009}. However, we follow the common practice and take $ B= 500$ in our calculation.

\subsection{Consistency of ePPR.aga estimator}

Two models are considered in the study of the theoretical properties of $m_k^{ens}$, both of which include the general regressions.
The first one is an extended additive model (XAM), which is an extension of the additive model:
\begin{equation}\label{Extendedmodel1}
  m(x)=\sum_{\tau=1}^L{m_\tau(x_{\C_\tau})},
\end{equation}
where each $m_\tau$ is a function of only $q$ variables $x_{\C_\tau}$ with $ Card(\C_\tau) = q $. It is worth noting that XAM can be any \textit{general function with $p$ variables} when $q=p$ and the corresponding ePPR estimator is exactly the one given in Section \ref{est. via AGA}.  The second model is the extended PPR model:
\begin{equation}\label{Pprmode1}
  m(x)=\sum_{\tau=1}^\infty {g_\tau(\theta_\tau^\top  x_{\C_\tau})},
\end{equation}
where $\theta_\tau\in\mathbb{R}^{q}$ and $ Card(\C_\tau )= q $.  Again, model \eqref{Pprmode1} includes the general PPR expansion \eqref{PPRidentity} when $ q=p$, and is called XPPR hereafter.
 When $ q = 1$, both models above degenerate to the additive model. For convenience of derivation, the theoretical results are only given for the truncated version of ${m}^{ens}_{k_n}$,
$
\hat{m}^{ens}_{k_n}=\frac{1}{B}\sum_{b=1}^{B} {\hat{m}_{k_n,b}},
$
where $\hat{m}_{k_n, b}= \max\{\min\{{m}_{k_n,b},t_n\},-t_n\}, b=1,\ldots, B$, and $t_n$ is the threshold, and $t_n\to\infty$, while in the calculation this truncation is not necessary. Asymptotic results about $\hat{m}^{ens}_{k_n }$ are concluded in Theorem \ref{EPPRcons} and \ref{EPPRresultppr} for the two models respectively.


\begin{theorem}\label{EPPRcons}
For model \eqref{Extendedmodel1}, assume
$
\E(e^{c\cdot Y^2})<\infty
$
for some $c>0$. Suppose $m$ is bounded and each $m_\tau\in L^2(\B_1^{q}), 1\leq \tau \leq L$.  If $t_n\asymp \ln n$ and $k_n,J_n\to\infty$ and $k_nJ_n\ln(k_nJ_n)/n\to 0$, then
$$
\E_{\mathcal{D}_{n},\Xi_{k_n}} \int|\hat{m}^{ens}_{k_n}(x)-m(x)|^2d\mu(x)\to 0.
$$
\end{theorem}

 Theorem \ref{EPPRcons} reveals that ePPR.aga works for XAM and is also consistent for any \textit{general function with $p$ variables} when  $q=p$ as presented in Theorem \ref{PPRpprmodelrate22}. Because XAM offers greater flexibility than the additive model, which is the main focus in the theoretical study of RF (see e.g. \cite{scornet2015consistency} and \cite{klusowski2022large}), ePPR.aga can be applied in a wider range of cases. The condition $k_nJ_n\ln(k_nJ_n)/n= o(1)$ in Theorem \ref{EPPRcons} guarantees that the variance of $\hat{m}^{ens}_{k_n}(x)$ is not large. A similar condition is also employed to restrict the size of ANN in the consistency study of ANN (see  e.g.  \cite{schmidt2020nonparametric} and \cite{kohler2021rate} ).


\begin{theorem}\label{EPPRresultppr}
For  model \eqref{Pprmode1}, assume
$
\E(e^{c\cdot Y^2})<\infty
$
for some $c>0$.  Suppose  $m\in\mathscr{M}_{s,q}$  with $\sum_{\iota=1}^\infty {\|m_\iota\|_\infty}<\infty$ and $\sum_{\iota=1}^\infty {TV(m_\iota^{(s)})}<\infty$.
When $t_n\asymp \ln n$, $k_n\asymp n^{(1+\frac{2s-1}{2q})/(2r)}$ and $J_n\asymp n^{(1+\frac{2s-1}{2q})^2/(4s^2)}$, we have
$$
\E_{\mathcal{D}_{n},\Xi_{k_n}} \int|\hat{m}^{ens}_{k_n}(x)-m(x)|^2d\mu(x) = O(
n^{-\alpha_*(s,q)}\cdot \ln^4 n),
$$
where  $\alpha_*(s,q)= (\frac{s+.5}{q}+1)/\left((\frac{s+.5}{q}+1)+ (1+\frac{1}{2s})(1-\frac{1}{2q})\right)$.
\end{theorem}

\begin{remark}

We have three remarks about the consistency rate $\alpha_*(s,q)$ of ePPR.aga.
Firstly, $\alpha_*(s,q)$ depends on $q$ but not $p$, if the PPR model is correct. Secondly, the proof of Theorem \ref{EPPRresultppr} is valid for Theorem \ref{PPRpprmodelrate22} if $ q = p$. Finally, as $s\to\infty$, $\alpha_*(s,q)\to 1$, which means that the consistency rate tends to the optimal  but is slightly reduced by a factor of $ \ln ^4 n $. On the other hand, it is known that RF is consistent under the model \eqref{Pprmode1} when $q=1$, i.e., the additive model \citep{scornet2015consistency}. However, to our best knowledge,  only \cite{klusowski2022large}  showed that the consistency rate of RF  is  $O_p\left((\ln n)^{-1}\right)$, which is much slower than the rate in Theorem \ref{EPPRresultppr}. More recently, \cite{chi2022asymptotic} proved that RF has  a polynomial consistency rate under the ``sufficient impurity decrease'' (SID) condition. However,  SID condition is only satisfied by some special cases, such as when $m(x),x\in\B_1^p$ is monotone along some directions. Therefore, ePPR.aga is more efficient than RF in the minimax sense, and its statistical consistency is guaranteed for more general regressions.
\end{remark}


\subsection{Stopping rule for the algorithms}
In practice, we still need a criterion to stop the algorithm at the right time  of the iterations  to balance the estimation bias and over-fitting. 
Next, we propose  a BIC criterion for this purpose. Consider the $b$-th run of the algorithm, denote by $m_{\tau,b}$ the output function of the inner loop of ePPR.aga. Define BIC  for each  individual estimator $m_{\tau,b},\tau=1,2,\ldots, $ by
\begin{equation}\label{BICeppr}
    \text{BIC}_{b}(\tau)= n^{-1} \sum_{i=1}^n (Y_i -m_{\tau,b}(X_i))^2+\nu\ln^5{n}\cdot\frac{\tau(q+J_n)}{n},
\end{equation}
where $\nu>0$. In the BIC above, $\tau (q+J_n) $ can be regarded as the number of parameters in the model. By including the factor $\ln^5{n}$ in \eqref{BICeppr}, the aim is to regulate the fluctuations induced by the sub-Gaussian error $Y-\E(Y|X)$, thereby guaranteeing proper  model selection and prediction. If the residual is bounded, the factor $\ln^5{n}$ can be reduced to $\ln{n}$, as demonstrated by a similar case of RF presented in
\cite[p.~1722]{scornet2015consistency}.  We choose the stopping time $k^*_{b}$ for the $b$-th run by
\begin{equation}\label{BICchooseK}
    k^*_{b}= \argmin_{1\leq \tau<\infty}\text{BIC}_{b}(\tau).
\end{equation}
In calculation, we select $k_b^*$ as the first $\tau\geq 1$ such that $BIC_b(\tau)<BIC_b(\tau+1)$. Finally, the ePPR estimator based on BIC criterion is given by
\begin{equation}\label{EppRBICestimator}
m^{ens}_{*}= \frac{1}{B}\sum_{b=1}^{B} {m_{k^*_{b},b}}.
\end{equation}
In \eqref{BICchooseK}, we need to minimize $BIC(\tau)$ over all $\tau\in\mathbb{Z}^+$.  In practice, we only need considering finite terms since $\nu\cdot{(k_b^*(q+J_n)\ln^5{n})}/{n} \leq \|\mathbb{Y}\|_n^2$. Again, to analyze asymptotic properties of estimator \eqref{EppRBICestimator}, the truncated version  $m^{ens}_{*}$:
$
\hat{m}^{ens}_{*}=\frac{1}{B}\sum_{b=1}^{B} {\hat{m}_{k^*_{b},b}}
$
is considered again where $\hat{m}_{k^*_{b},b}= \max\{\min\{m_{k^*_{b},b},t_n\},-t_n\}$ and $t_n>0$ is the threshold.

\begin{theorem}\label{EPPRcons2a}
For model \eqref{Extendedmodel1}, assume
$
\E(e^{c\cdot Y^2})<\infty,
$
where the constant $c>0$. Suppose $m$ is bounded and each $m_\tau\in L^2(\B_1^{q}), 1\leq \tau \leq L$.  There is $\nu_0>0$ such that if $\nu\ge \nu_0$ and $t_n\asymp \ln n$, $J_n=o\left(\frac{\ln^6 n}{n}\right)$ and $J_n\to\infty$, then
$$
 \int|\hat{m}^{ens}_{*}(x)-m(x)|^2d\mu(x)=o_p(1).
$$
\end{theorem}

Note that Theorem \ref{EPPRcons}   requires that $k_nJ_n\ln(k_nJ_n)/n=o(1)$, while a stronger condition that factor $\ln^6{n}$ in $J_n=o\left(\frac{\ln^6 n}{n}\right)$ is necessary in Theorem \ref{EPPRcons2a}. This is because $\hat{m}^{ens}_{*}$ is obtained by imposing a penalty shown in \eqref{BICeppr}, which depends on $\nu$.

\begin{theorem}\label{EPPRpprmodelrate2b}
For  model \eqref{Pprmode1}, assume
$
\E(e^{c\cdot Y^2})<\infty
$
for some  $c>0$.  Suppose  $m\in\mathscr{M}_{s,q}$  with $\sum_{\iota=1}^\infty {\|m_\iota\|_\infty}<\infty$ and $\sum_{\iota=1}^\infty {TV(m_\iota^{(s)})}<\infty$.  There exists $\nu_0>0$ such that if $\nu\ge \nu_0$ and $t_n\asymp \ln n$ and $J_n\asymp n^{(1+\frac{2r-1}{2q})^2/(4s^2)}$,  we have
$$
 \int|\hat{m}^{ens}_{*}(x)-m(x)|^2d\mu(x) = O_p(
n^{-\alpha_*(s,q)}\cdot \ln^5 n),
$$
where rate $\alpha_*(s,q)$ is given in Theorem \ref{EPPRresultppr}.
\end{theorem}

Theorem \ref{EPPRpprmodelrate2b} reveals that $\hat{m}^{ens}_{*}$ selected by BIC criterion has a similar consistency rate as $\hat{m}^{ens}_{k_n}$ in Theorem \ref{EPPRresultppr}.
 The difference is that the rate of $\hat{m}^{ens}_{*}$ depends on the penalty parameter $\nu$  defined in \eqref{BICeppr}.
Although the choice of $\nu_0$ in the theorem depends on various factors, our experience suggests that setting $\nu = 2$ yields satisfactory results. Additionally, it is worth noting that the penalty parameter is actually $\nu \log^5 n$, which is slightly higher than the conventional BIC. This adjustment is supported by our numerical study as well.

\section{Numerical Performance in Real Data} \label{secReal}

In this section, we use 36 real data sets with continuous responses  and 36 data sets with binary categorical response (0 and 1) to demonstrate the performance of ePPR for regression prediction and classification respectively.
These data sets can be found on one of the following websites:
(A) \url{https://archive.ics.uci.edu/ml/datasets};
(B) \url{https://www.kaggle.com};
(C) \url{http://lib.stat.cmu.edu/datasets};
(D) \url{https://data.world};
(E) \url{https://www.openml.org}.

\begin{table}[p!]\small
  \renewcommand\arraystretch{1.15}
  \caption{Regression: average RPE based on 100 random partitions of each data set into training and test sets}\label{Table1}
\centerline{\hspace{-0.1cm}
\begin{tabular}{l@{\hspace{0.5\tabcolsep}}r@{\hspace{0.5\tabcolsep}}r@{\hspace{0.5\tabcolsep}}|
r@{\hspace{0.2\tabcolsep}}r@{\hspace{0.6\tabcolsep}}r@{\hspace{0.6\tabcolsep}}
r@{\hspace{0.6\tabcolsep}}r@{\hspace{0.6\tabcolsep}}r@{\hspace{0.6\tabcolsep}}
r@{\hspace{0.6\tabcolsep}}r}
\hline
  Data sets (and sources)	&	$N$	&	$p$	& \multicolumn{1}{c}{LM}	&	\multicolumn{1}{c}{PPR}	&	SVM	&	\multicolumn{1}{c}{RF}	&	GRF	&	RLT	&	XGB	&	{ ePPR}	\\\hline
 { Concrete compressive strength (A)}	&	1030	&	8	&	0.270	&	0.142	&	0.167	&	0.120	&	0.169	&	0.204	&	0.083	&	\textbf{\small  0.079}	\\
  Wine quality (B)	&	1143	&	11	&	0.649	&	0.675	&	0.612	&	\textbf{\small  0.559}	&	0.626	&	0.580	&	0.644	&	0.607	\\
  Boston house price (C)	&	506	&	13	&	0.252	&	0.227	&	0.184	&	0.143	&	0.206	&	0.175	&	0.154	&	\textbf{0.125}	\\
  Wild blueberry yield (B)	&	777	&	13	&	5.424	&	18.64	&	0.043	&	0.074	&	0.121	&	0.120	&	0.052	&	\textbf{\small  0.028}	\\
  Body fat (B)	&	252	&	14	&	\textbf{\small  0.028}	&	0.061	&	0.126	&	0.101	&	0.058	&	0.132	&	0.052	&	0.038	\\
  Baseball player salary (C)	&	263	&	16	&	0.483	&	0.605	&	0.302	&	0.249	&	0.304	&	0.270	&	0.315	&	\textbf{\small  0.247}	\\
  Paris housing price (B)	&	10000	&	16	&	0.001	&	\textbf{\small  0.000}	&	0.026	&	0.019	&	\textbf{\small  0.000}	&	0.054	&	\textbf{\small  0.000}	&	\textbf{\small  0.000}	\\
  Insurance churn participants (B)	&	33908	&	16	&	0.794	&	0.921	&	0.886	&	\textbf{\small  0.715}	&	0.778	&	0.743	&	0.857	&	0.758	\\
  House sales in King County, USA (B)	&	21613	&	18	&	0.228	&	0.217	&	0.226	&	0.161	&	0.199	&	0.195	&	0.157	&	\textbf{\small  0.134}	\\
  Parkinsons (A)	&	5875	&	19	&	0.823	&	0.759	&	0.671	&	0.352	&	0.387	&	0.524	&	\textbf{\small  0.125}	&	0.224	\\
  Life Expentency (WHO) (B)	&	1649	&	20	&	0.136	&	0.086	&	0.095	&	\textbf{\small  0.049}	&	0.074	&	0.064	&	0.057	&	0.056	\\
  Bar crawl (B)	&	7590	&	21	&	0.231	&	0.225	&	0.193	&	0.170	&	0.235	&	0.192	&	0.172	&	\textbf{\small  0.144}	\\
  Cancer death rate (D)	&	2332	&	32	&	0.522	&	0.614	&	0.536	&	0.502	&	0.556	&	0.524	&	0.524	&	\textbf{\small  0.187}	\\
  Wave energy converters-Adelaide (A)	&	71998	&	32	&	0.841	&	0.614	&	\textbf{\small  0.203}	&	0.280	&	0.453	&	0.302	&	0.326	&	0.237	\\
  Wave energy converters-Perth (A)	&	72000	&	32	&	0.885	&	0.708	&	\textbf{\small  0.258}	&	0.345	&	0.561	&	0.375	&	0.395	&	0.290	\\
  Sidney house price (B)	&	30000	&	37	&	0.638	&	1.000	&	0.631	&	\textbf{\small  0.555}	&	0.618	&	0.559	&	0.649	&	0.564	\\
  Russian house price (B)	&	21329	&	50	&	\textbf{\small  0.647}	&	1.238	&	0.678	&	0.660	&	0.671	&	0.653	&	0.789	&	0.653	\\
  Facebook comment volume (A)	&	18370	&	52	&	0.727	&	2.225	&	0.808	&	0.567	&	0.630	&	0.625	&	0.648	&	\textbf{\small  0.536}	\\
  Gold price prediction (B)	&	1718	&	74	&	0.018	&	0.011	&	0.009	&	\textbf{\small  0.005}	&	0.010	&	\textbf{\small  0.005}	&	0.008	&	0.006	\\
  Baseball player statistics (B)	&	4535	&	74	&	0.001	&	\textbf{\small  0.000}	&	0.059	&	0.019	&	0.027	&	0.033	&	0.002	&	0.001	\\
  CNNpred (A)	&	1441	&	76	&	0.001	&	\textbf{\small  0.000}	&	0.056	&	0.001	&	0.003	&	0.002	&	0.002	&	\textbf{\small  0.000}	\\
  Superconductivity (A)	&	21263	&	81	&	0.285	&	0.257	&	0.226	&	\textbf{\small  0.161}	&	0.221	&	0.169	&	0.180	&	0.184	\\
  Warsaw flat rent price (B)	&	3472	&	83	&	0.333	&	0.461	&	0.488	&	0.308	&	0.370	&	0.343	&	0.331	&	\textbf{\small  0.277}	\\
  Year prediction MSD (A)	&	50000	&	90	&	0.785	&	1.635	&	0.796	&	0.853	&	0.906	&	0.850	&	1.001	&	\textbf{\small  0.754}	\\
  Buzz in social media (A)	&	28179	&	96	&	0.200	&	0.460	&	0.655	&	0.103	&	0.204	&	0.132	&	0.132	&	\textbf{\small  0.103}	\\
  Communities and crime (A)	&	1994	&	101	&	0.348	&	0.746	&	0.372	&	0.347	&	0.374	&	0.344	&	0.397	&	\textbf{\small  0.326}	\\
  Residential building-Sales (A)	&	372	&	103	&	0.027	&	0.043	&	0.124	&	0.046	&	0.081	&	0.069	&	0.019	&	\textbf{\small  0.009}	\\
  Residential building-Cost (A)	&	372	&	103	&	0.029	&	0.232	&	0.067	&	0.064	&	0.101	&	0.081	&	0.050	&	\textbf{\small  0.026}	\\
  Geographical original-Latitude (A)	&	1059	&	116	&	0.827	&	1.471	&	0.765	&	\textbf{\small  0.737}	&	0.831	&	0.759	&	0.823	&	0.748	\\
  Geographical original-Longitude (A)	&	1059	&	116	&	0.739	&	1.200	&	\textbf{\small  0.654}	&	0.673	&	0.791	&	0.708	&	0.733	&	0.659	\\
  Blog feedback (A)	&	13063	&	249	&	0.645	&	1.136	&	0.742	&	\textbf{\small  0.556}	&	0.613	&	0.588	&	0.650	&	0.559	\\
  Credit score (B)	&	80000	&	259	&	0.208	&	0.497	&	0.179	&	0.128	&	0.178	&	0.126	&	0.131	&	\textbf{\small  0.111}	\\
  CT slices (A)	&	53500	&	380	&	0.206	&	0.214	&	0.087	&	0.080	&	0.173	&	\textbf{\small  0.075}	&	0.094	&	0.088	\\
  UJIndoor-Longitude (A)	&	19937	&	465	&	0.099	&	0.143	&	0.087	&	\textbf{\small  0.020}	&	0.043	&	0.030	&	0.037	&	0.030	\\
  UJIndoor-Latitude (A)	&	19937	&	465	&	0.131	&	0.341	&	0.121	&	\textbf{\small  0.033}	&	0.083	&	0.066	&	0.052	&	0.054	\\
  Regression with categorical data (B)	&	169486	&	1032	&	0.621	&	6.957	&	0.701	&	0.576	&	0.659	&	0.583	&	0.635	&	\textbf{\small  0.556}	\\\hdashline
  \multicolumn{3}{c|}{Average of RPE across all data sets}&0.530	&	1.244	&	0.357	&	0.287	&	0.342	&	0.313	&	0.313	&	{0.261}	\\
  \multicolumn{3}{c|}{\textit{no. of bests}$^*$  }& 2 &1.75 &3 &9.5 &0.25 &1.5 &1.25 & {16.75}\\
    \hline
    \multicolumn{11}{l}{ * \footnotesize   the number of data sets for which the method performs the best among all the methods;  if $Q$ methods have the }\\
     \multicolumn{11}{l}{\footnotesize same  best performance for one data set, count the no. as $1/Q$ for each of the methods.}
  \end{tabular}
 }
\end{table}

\begin{table}[h!]\small
  \renewcommand\arraystretch{1.15}
   \centering
   \caption{Classification: average MR(\%) based on 100 random partitions of each data set into training and test sets}
   \label{Table2}
\centerline{  \hspace{-0.4cm}
 \begin{tabular}{l@{\hspace{0.5\tabcolsep}}r@{\hspace{0.5\tabcolsep}}r@{\hspace{0.5\tabcolsep}}|
r@{\hspace{0.5\tabcolsep}}r@{\hspace{0.5\tabcolsep}}r@{\hspace{0.5\tabcolsep}}
r@{\hspace{0.5\tabcolsep}}r@{\hspace{0.5\tabcolsep}}r@{\hspace{0.5\tabcolsep}}r@{\hspace{0.5\tabcolsep}}r
@{\hspace{0.5\tabcolsep}}r@{\hspace{0.5\tabcolsep}}r}
  \hline
  Data sets (and sources)	&	$N$	&	$p$	&	GLM	&	PPR	&	SVM	&	RF	&	GRF	&	RLT	&	XGB	&	Ada	&	ANN	&	{ePPR}	\\\hline
  Indian liver patient (A)	&	579	&	10	&	29.52	&	31.56	&	29.16	&	29.67	&	\textbf{\small  28.82}	&	29.25	&	32.02	&	30.27	&	29.22	&	29.31	\\
  MAGIC Gamma telescope (A)	&	19020	&	10	&	20.58	&	15.85	&	15.62	&	14.68	&	16.42	&	14.73	&	17.90	&	14.89	&	36.97	&	\textbf{\small  13.73}	\\
  Heart disease (A)	&	270	&	13	&	\textbf{\small  16.41}	&	21.53	&	17.29	&	17.47	&	17.99	&	17.06	&	22.81	&	17.92	&	44.11	&	18.20	\\
  Australian credit approval (A)	&	690	&	14	&	13.75	&	14.76	&	14.57	&	\textbf{\small  13.12}	&	14.23	&	13.50	&	15.68	&	13.27	&	14.96	&	13.44	\\
  EEG eye state (A)	&	14980	&	14	&	37.40	&	24.18	&	25.39	&	18.04	&	24.74	&	19.30	&	22.57	&	20.87	&	44.84	&	\textbf{\small  11.05}	\\
  seismic-bumps (A)	&	2584	&	15	&	6.62	&	7.47	&	6.59	&	6.78	&	6.56	&	6.64	&	9.08	&	6.85	&	\textbf{\small  6.50}	&	7.34	\\
  Autism screening adult (A)	&	702	&	16	&	\textbf{\small  0.00}	&	0.04	&	2.36	&	\textbf{\small  0.00}	&	\textbf{\small  0.00}	&	0.28	&	\textbf{\small  0.00}	&	\textbf{\small  0.00}	&	15.90	&	0.01	\\
  Bank marketing (A)	&	45211	&	16	&	11.18	&	11.92	&	11.09	&	\textbf{\small  10.61}	&	11.50	&	10.92	&	12.27	&	10.91	&	11.94	&	10.91	\\
  Diabetic retinopathy debrecen (A)	&	1151	&	19	&	26.05	&	26.80	&	30.78	&	31.74	&	33.51	&	32.38	&	34.87	&	31.41	&	48.39	&	\textbf{\small  25.12}	\\
  Performance prediction (B)	&	1340	&	19	&	\textbf{\small  29.33}	&	31.43	&	29.50	&	31.27	&	29.94	&	30.45	&	36.26	&	31.14	&	40.07	&	29.78	\\
  Climate model simulation crash (E)	&	540	&	20	&	8.99	&	12.91	&	\textbf{\small  8.58}	&	8.71	&	\textbf{\small  8.58}	&	8.59	&	11.21	&	9.12	&	8.72	&	9.27	\\
  Default of credit card clients (A)	&	30000	&	24	&	19.68	&	20.43	&	19.15	&	19.15	&	\textbf{\small  18.30}	&	18.97	&	24.13	&	19.22	&	22.12	&	19.24	\\
  Parkinson multiple sound (B)	&	1208	&	26	&	33.06	&	33.85	&	29.60	&	27.61	&	30.39	&	28.37	&	30.83	&	28.49	&	43.25	&	\textbf{\small  26.02}	\\
  Pistachio (B)	&	2148	&	28	&	7.50	&	8.02	&	8.19	&	10.80	&	12.63	&	10.91	&	11.69	&	9.34	&	44.36	&	\textbf{\small  6.57}	\\
  Breast cancer (B)	&	569	&	30	&	2.95	&	5.46	&	2.67	&	4.24	&	5.97	&	4.02	&	5.52	&	3.85	&	37.95	&	\textbf{\small  2.41}	\\
  Phishing websites (A)	&	11055	&	30	&	7.80	&	8.73	&	6.86	&	\textbf{\small  5.67}	&	8.01	&	6.61	&	6.81	&	6.25	&	7.74	&	6.29	\\
  Ionosphere (A)	&	351	&	33	&	12.65	&	15.75	&	6.14	&	7.21	&	7.65	&	7.25	&	10.79	&	7.52	&	12.65	&	\textbf{\small  5.61}	\\
  QSAR biodegradation (A)	&	1055	&	41	&	14.13	&	16.16	&	13.02	&	13.38	&	16.08	&	13.18	&	15.85	&	13.88	&	30.23	&	\textbf{\small  12.55}	\\
  Spambase (A)	&	4601	&	57	&	7.04	&	9.19	&	8.47	&	6.16	&	8.30	&	6.29	&	8.35	&	6.32	&	41.62	&	\textbf{\small  5.87}	\\
  Mice protein expression (A)	&	1047	&	70	&	2.94	&	3.89	&	1.00	&	1.35	&	7.47	&	1.28	&	6.89	&	3.10	&	6.70	&	\textbf{\small  0.54}	\\
  Ozone level detection (A)	&	1847	&	72	&	6.66	&	11.72	&	7.00	&	6.44	&	7.02	&	6.70	&	7.45	&	6.37	&	7.06	&	\textbf{\small  6.35}	\\
  Insurance company benchmark (A)	&	5822	&	85	&	\textbf{\small  5.99}	&	15.22	&	\textbf{\small  5.99}	&	6.54	&	\textbf{\small  5.99}	&	6.20	&	9.03	&	6.33	&	6.16	&	6.77	\\
  Company bankruptcy (B)	&	6819	&	94	&	3.21	&	5.34	&	3.22	&	\textbf{\small  3.13}	&	3.22	&	\textbf{\small  3.13}	&	3.81	&	3.15	&	3.25	&	3.39	\\
  Hill-valley (A)	&	1212	&	100	&	50.48	&	4.32	&	50.00	&	40.76	&	46.00	&	41.93	&	38.13	&	43.15	&	51.93	&	\textbf{\small  1.58}	\\
  Musk (A)	&	6598	&	166	&	7.22	&	10.24	&	5.68	&	5.63	&	9.63	&	6.21	&	7.06	&	6.15	&	15.45	&	\textbf{\small  4.34}	\\
  ECG heartbeat categorization (B)	&	14550	&	186	&	19.18	&	25.45	&	18.06	&	8.88	&	14.86	&	8.88	&	13.36	&	10.33	&	13.26	&	\textbf{\small  7.30}	\\
  Arrhythmia (A)	&	420	&	192	&	24.79	&	37.20	&	23.64	&	19.67	&	23.05	&	21.86	&	24.89	&	\textbf{\small  18.76}	&	45.29	&	20.19	\\
  Financial indicators of US stocks (B)	&	986	&	216	&	0.62	&	12.93	&	13.94	&	0.94	&	1.33	&	4.74	&	\textbf{\small  0.05}	&	\textbf{\small  0.05}	&	23.53	&	0.31	\\
  Madelon (A)	&	2000	&	500	&	38.91	&	48.10	&	42.46	&	34.52	&	37.54	&	32.98	&	32.78	&	\textbf{\small  32.46}	&	49.71	&	33.44	\\
  Human activity recognition (A)	&	2633	&	561	&	\textbf{\small  0.00}	&	\textbf{\small  0.00}	&	0.10	&	0.08	&	0.25	&	0.16	&	0.01	&	0.01	&	0.01	&	\textbf{\small  0.00}	\\
  Gina agnostic (E)	&	3468	&	970	&	14.38	&	38.87	&	11.47	&	9.25	&	13.36	&	9.20	&	13.80	&	\textbf{\small  9.18}	&	45.33	&	9.93	\\
  QSAR androgen receptor (A)	&	1687	&	1024	&	10.36	&	35.13	&	10.09	&	9.71	&	11.38	&	\textbf{\small  9.55}	&	11.69	&	9.71	&	11.38	&	9.59	\\
  QSAR oral toxicity (A)	&	8992	&	1024	&	8.10	&	36.98	&	7.71	&	\textbf{\small  7.26}	&	8.18	&	7.31	&	9.08	&	7.42	&	8.33	&	7.33	\\
  Predicting a biological response (B)	&	3751	&	1776	&	25.20	&	34.10	&	25.98	&	23.15	&	25.17	&	\textbf{\small  22.45}	&	27.55	&	23.60	&	33.06	&	23.62	\\
  GISETTE (A)	&	6000	&	4955	&	4.84	&	16.36	&	5.29	&	5.16	&	8.34	&	4.77	&	7.72	&	4.69	&	39.35	&	\textbf{\small  4.30}	\\
  Arcene (A)	&	200	&	9961	&	28.69	&	44.27	&	25.93	&	19.76	&	27.93	&	22.03	&	28.10	&	21.84	&	38.06	&	\textbf{\small  13.69}	\\\hdashline
  \multicolumn{3}{c|}{average of MR(\%) across all data sets}&15.45	&	19.34	&	15.07	&	13.29	&	15.29	&	13.56	&	15.83	&	13.55	&	26.09	&	11.26	\\
  \multicolumn{3}{c|}{\textit{no. of bests}}&2.87 &0.33 &0.83 &4.7 &3.03 &2.5 &0.7 &3.7 &1 &16.33\\
   \hline
  \end{tabular}
  }
\end{table}

\begin{table}[h!]\small
  \caption{Average RPE (or MR \%) for each method across the 36 variable-augmented regression data sets (or 36 variable-augmented classification data sets) with p = 10000 and the number of data sets for which the method has the best performance (\textit{no. of bests}) among all the methods}
  \label{Table3}
  \centering
  \begin{tabular}{cc|c@{\hspace{1\tabcolsep}}c@{\hspace{1\tabcolsep}}c@{\hspace{1\tabcolsep}}c@{\hspace{1\tabcolsep}}c@{\hspace{1\tabcolsep}}c@{\hspace{1\tabcolsep}}c}\hline
   &  &LM&SVM&RF&GRF&RLT&XGB&{\footnotesize ePPR}\\\hline
the 36  regression   & {average RPE}&0.469&0.860&0.394&0.748&0.446&0.407&0.348\\ \cline{2-9}
 data sets  & \textit{no. of bests}& 1 &0 &6 &0 &1 &2 &26\\
    \hline
the 36 classification  &    {average MR(\%)}&17.95&28.12&22.97&25.45&17.23&20.17&15.90\\     \cline{2-9} data sets & \textit{no. of bests}& 5.06&0.76&0.76&0.76&6.56&4.5&17.56\\
 \hline
  \end{tabular}
\end{table}


For regression analysis, we compare the performance of ePPR with several other methods, including the original PPR, support vector machine (SVM), random forest (RF), generalized random forest (GRF), reinforcement learning trees (RLT), and extreme gradient boosting (XGB).
In the context of classification, we evaluate ePPR against PPR, SVM, RF, GRF, RLT, XGB, adaptive boosting (Ada), and artificial neural network (ANN). To apply PPR, GRF, and ePPR to classification tasks, which were originally developed for regression, we treat the class labels (0 and 1) as continuous values. We then utilize these methods directly to predict the response variable, denoted as $\hat{y}$. Finally, we classify the data based on the threshold $0.5$, using the indicator function $\mathbb{I}(\hat{y} > 0.5)$, where $\mathbb{I}(\cdot)$ represents the indicator function. It is worth noting that extending these methods to multiclass classification can also be made.

The following functions and packages in R are used for the calculation of the methods: \texttt{glmnet} \citep{glmnetpackage} is used for Lasso, including LM for linear regression and GLM (logistic model) for classification, \texttt{svm} function in package \texttt{e1071} \citep{e1071package} for SVM,  \texttt{randomForest} \citep{randomForestpackage} for RF,  \texttt{regression\_forest} in package \texttt{grf} \citep{athey2019generalized} for GRF,  \texttt{RLT} \citep{Zhu2015} for RLT,  package \texttt{xgboost} \citep{xgboostpackage} for XGB, \texttt{ada} \citep{adapackage} for ada, and \texttt{mlp} function in package \texttt{RSNNS} \citep{RSNNSpackage} for ANN. For all the R functions  and  packages, their default values of tuning parameters are used.  However, when it comes to ANN, there is no default values, i.e. number of hidden layers or nodes within each layer.  We use 5-fold CV to select the size from candidate sizes of an equal number of neurons $w$ in all hidden layers $L$ with $ L = 1, ..., 10 $ and $  w= 5, 6, ..., 50 $. For ePPR, $B=500$ is used.

For each data set, we randomly partition it into training set and test set. The training set consists of $n=\min(\lfloor 2N/3\rfloor,1000)$ randomly selected observations, where $N$ is the number of observations in the original data sets, and the remaining observations form the test set. For regression, the relative prediction error, defined as
$$RPE=\sum_{i\in \text{test set}}(\hat{y}_i-y_i)^2/\sum_{i\in \text{test set}}(\bar{y}_{\text{train}}-y_i)^2,$$
where $\bar{y}_{\text{train}}$ is a naive prediction based on the average of $y$ in the training sets, is used to measure the performance of a method. For classification, the misclassification rate, defined as
$$MR=\sum_{i\in \text{test set}} 1(1(\hat{y}_i>0.5) \neq y_i) /(N-n),$$
is used to measure the performance. For each data set, the random partition is repeated 100 times,  and averages of the RPEs or MRs are calculated to compare different methods. The calculation results are listed in {Table \ref{Table1}} and {Table \ref{Table2}}. The smallest RPE or MR for each data set is highlighted in \textbf{bold} font.

By comparing the prediction errors, RPEs or MRs, of all the methods, we have the following observations. The classical methods SVM, RF and XGB have reasonably good performance on nearly all the data sets in comparison with LM and PPR. In particular, RF has been commonly recognized as the most robust and efficient learning method, especially for tabular data; our extensive calculation also supports the recognition. The recent variants of RF, such as GRF and RLT, also show their effectiveness in the prediction.
Although PPR has good performance on some data sets, it is unstable and has much bigger prediction errors than the other methods for most of the data sets.  This could be the reason why PPR is not popularly received in practice.  The proposed ePPR is quite stable and attains the smallest RPE and MR in most data sets as indicated in {Table \ref{Table1}} and {Table \ref{Table2}}. The competency of ePPR is also verified by the fact that it has the smallest average of RPEs (or MSs) across all the 36 data sets among all the methods. The numbers of data sets for which a method performs the best among all the competitors, denoted by \textit{no. of bests}, also suggest that ePPR outperforms the other methods for most of the data sets. The reason for the better performance of ePPR over RF might be the fact that ePPR is smooth, while RF is a step function \cite{wager2015adaptive}.

It is also noticed that RF has good performance for high dimensional data \citep{caruana2006empirical, caruana2008empirical}. To see the performance of ePPR in high dimension, we add randomly generated i.i.d. normal variables, called augmented variables, to the original data  so that altogether there are $p=10000$ predictors,  and then compare all the methods. Since PPR, Ada and ANN require huge memory with big $ p $, they are not included for comparison in {Table \ref{Table3}}. For the other methods, the overall performances across all the data sets are listed in {Table \ref{Table3}}. We can see that ePPR still has the best performance in most of the data sets for both regression and classification. Although all methods suffer from the high dimension, ePPR is less influenced: compared to RF, ePPR has about 9\% improvement in prediction and 15\% improvement in classification in {Table \ref{Table1}} and {Table \ref{Table2}} respectively, while ePPR is about 12\% or 30\% respectively better than RF in {Table \ref{Table3}}, indicating that the relative improvement increases  with the dimension of the data.  The number of data sets for which ePPR outperforms the competitors,  \textit{no. of bests}, also increases with the dimension of the data.


\def\I{\mathcal{I}}
\def\MSE{\mathbf{MISE}}

\section{Proofs} \label{SecProof}
 Let $c$ be a positive constant which may change from line to line. Function $c(.)$ is positive and only depends on its arguments. 
For any estimator  $\hat m(x) $ of $ m(x) $, denote the integrated squared error (ISE) by
$$  {{\cal I}}(\hat m(x)) =  \int|\hat m(x) -m(x)|^2d\mu(x),  $$
and mean integrated squared error (MISE) by  $ \MSE(\hat m(x)) = \E_{\mathcal{D}_{n}} {{\cal I}}(\hat{m}(x)) $. For the ensemble estimator, we also need to take the average over the random subsets space $ \Xi_{k_n} $, the corresponding MISE is also denoted by
$$
\MSE (\hat{m}(x)) = \E_{\mathcal{D}_{n},\Xi_{k_n}} {{\cal I}}(\hat m(x)).
$$
Let $x_1,\ldots, x_n \in\mathbb{R}^p$ and set $x_1^n=\{x_1,\ldots, x_n\}$. Let $\mathcal{G}$ be a class of functions $g:\mathbb{R}^p\to\mathbb{R}$. An $L_\ell$-$\epsilon$-cover of $\mathcal{G}$ on $x_1^n$ is defined in \cite{gyorfi2002distribution}. Denote by $\mathcal{N}_\ell(\epsilon,\mathcal{G},x_1^n)$ the covering number. In this paper, $ \ell = 1 $ is used. We further simplify the covering number as  $\mathcal{N}(\epsilon,\mathcal{G},x_1^n)$; see \cite{gyorfi2002distribution} or Definition (S.1) 
for
more details in  \SM.





\vspace{0.3cm}

\textit{Proof of Theorem \ref{Aga population version}.}
First, recall that the (dyadic) entropy numbers of any function class $\mathscr{F}\subseteq \mathscr{H}$ is defined by
\begin{equation}\label{Entropy Number General}
   \epsilon_k(\mathscr{F})_{\mathscr{H}}:=\inf \{ \epsilon: \mathscr{F}\ \text{can be covered by }\ 2^k\ \text{balls of radius}\ \epsilon\},
\end{equation}
which measures the size of $\mathscr{F}$ or how much information  $\mathscr{F}$  carries on. For each index $1\le j\le L$, where $L\in\mathbb{Z}^+$, choose any
 $s_j(\theta_j^Tx)=\sum_{\iota=1}^{J_n}{a_{j,\iota}\cdot \sigma_s(\theta_j^Tx+b_{j,\iota})}\in Dic^s_{n}$ satisfying $\sum_{\iota=1}^{J_n} |a_{j,\iota}| \leq 1$ and $b_{j,\iota}\in [-2,2]$.  Then, for any weights $w_j\in\mathbb{R}, j=1,\ldots, L$ with  $\sum_{j=1}^{J_n}{|w_j|}\leq 1$, we have
$$
 \sum_{j=1}^L{w_j\cdot s_j(\theta_j^Tx)}=\sum_{j=1}^L\sum_{\iota=1}^{J_n}{w_ja_{j,\iota}\cdot \sigma_s(\theta_j^Tx+b_{j,\iota})}
$$
satisfying
$
\sum_{j=1}^L\sum_{\iota=1}^{J_n}{|w_ja_{j,\iota}|}\leq 1.
$
Therefore,  we have
\begin{equation}\label{dgjhasdgba..}
\sum_{j=1}^L{w_j\cdot s_j(\theta_j^Tx)}\in B_1(\mathbb{P}_{s,2}^p),
\end{equation}
where the dictionary $\mathbb{P}_{s,2}^p$ is a class of neural networks:
$
  \mathbb{P}_{s,2}^p:= \{ \sigma_s(\theta^T x+b): \theta\in\Theta^p, b\in[-2,2], x\in\B_1^p\}.
$
Thus, \eqref{dgjhasdgba..} implies $$B_1(Dic_n^s)\subseteq B_1(\mathbb{P}_{s,2}^p)$$ for any $n\in\mathbb{Z}^+$. By equation $(3.51)$ in \cite{siegel2021sharp}, we have
$$
  \epsilon_k(B_1(Dic^s_{n}))_{L^2(\mathbb{B}_1^p,X)}\leq \epsilon_k(B_1(\mathbb{P}_{s,2}^p))_{L^2(\mathbb{B}_1^p)}\leq k^{-\frac{1}{2}- \frac{2s+1}{2p}}.
$$
As a consequence, the proof is completed by using Theorem 1 in \cite{siegel2022optimal}.   \hfill\(\Box\)





\vspace{0.3cm}

\textit{Proof of Theorem \ref{Lower bound of Aga}.}
Since $J_n$ is fixed when $n$ is given, without loss of generality, we only consider the case $J_n= 1$ in this proof.  Note that there exits some constant $M>0$ such that  $\sup_{x\in\B_1^p}{d_n(x)}\leq M$ for any $d_n\in Dic^s_{n}$, where $M$ does not depend on $n\in\mathbb{Z}^+$.  Based on this  fact, we can construct a better approximation based on AGA  by truncating the output function $m^{aga}_k$.  Specifically, consider the function as follows:
$$
T_{t_k} m^{aga}_k(x):= \max\{\min\{m^{aga}_k(x),t_n\},-t_n\},
$$
where we set $t_k=O(\log k)$.  Then, we have
$$
  \sup_{m\in B_1(Dic^s_{n})}{\|m-T_{t_n} m^{aga}_k(x)\|_{L^2(\mathbb{B}_1^p,X)}}\leq  \sup_{m\in B_1(Dic^s_{n})}{\|m- m^{aga}_k(x)\|_{L^2(\mathbb{B}_1^p,X)}}.
$$
Next, we show
\begin{equation}\label{optimal aga rate population bb}
 \sup_{m\in B_1(Dic^s_{n})}\sup_{k\in\mathbb{Z}^+}k^{-\frac{1}{2}-\alpha'(s,p)}\inf_{m_{k}\in \Sigma^{t_n}_{k}(Dic^s_{n})}{\|m(X)-m_{k}(X)\|_{L^2(\mathbb{B}_1^p,X)}}=\infty,
\end{equation}
where $\Sigma^{t_k}_{k}(Dic^s_{n}):= \{T_{t_k} g_n: g_n\in \Sigma_{k}(Dic^s_{n})\}$ is the truncated  version of $\Sigma_{k}(Dic^s_{n})$,  a linear space spanned by $Dic^s_{n}$.  Therefore,  rate $\alpha(s,p)$ in \eqref{HJKSBKJBKNKJNKNqs} can not be improved further and  AGA  achieves this optimal function approximation rate.

We use the proof by contradiction to show \eqref{optimal aga rate population bb}.  Suppose \eqref{optimal aga rate population bb} is false,  then there exists some constant $c_1>0$ such that
\begin{equation}\label{odchbkasjdhkadhqqzz.}
    \sup_{m\in B_1(Dic^s_{n})}\inf_{m_{k}\in \Sigma^{t_k}_{k}(Dic^s_{n})}{\|m_n-m_{k}\|_{L^2(\mathbb{B}_1^p,X)}} \leq c_1\cdot  k^{-\frac{1}{2}-\alpha'(s,p)}
\end{equation}
for any $k\in \mathbb{Z}^+$.  Recall that $\mathbb{P}_{s,2}^p$  appearing in the proof of  Theorem \ref{Aga population version} is
$$
  \mathbb{P}_{s,2}^p:= \{ \sigma_s(\theta^T x+b): \theta\in\Theta^p, b\in[-2,2], x\in\B_1^p\}
$$
and its convex symmetric hull is
$$
    Conv(\mathbb{P}_{s,1}^p):= \left\{  \sum_{\iota=1}^L{\alpha_\iota\cdot
    \sigma_s(\theta_\iota^T x+b_\iota)}:  \sum_{\iota=1}^L{|\alpha_\iota|}=1, \sigma_s(\theta_\iota^T x+b_\iota)\in \mathbb{P}_{s,2}^p, L\in\mathbb{Z}^+ \right\}.
$$
For any $f=\sum_{\iota=1}^L{\alpha_\iota\cdot \sigma_s(\theta_\iota^T x+b_\iota)}\in Conv(\mathbb{P}_{s,2}^p)$,  we know each $\sigma_s(\theta_\iota^T x+b_\iota)\in Dic^s_{n}$.  Considering linear combinations of $k$ elements in $\mathbb{P}_{s,1}^p$, we have a class of shallow neural networks:
$$
\mathcal{N}\mathcal{N}_{k}:= \left\{ \sum_{\iota=1}^k{\beta_\iota\cdot \sigma_s(\theta_\iota^T x+b_\iota)}: \beta_\iota\in\mathbb{R}, b_\iota\in [-1,1], \theta\in \Theta^p\right\}.
$$
Meanwhile, its truncated version is also required in the proof,
$$
\mathcal{N}\mathcal{N}_{k}^{t_k}:=\{T_{t_k}g: g\in \mathcal{N}\mathcal{N}_{k}\},
$$
where $T_{t_k}$ is the truncated operator defined above. From \eqref{odchbkasjdhkadhqqzz.},  there is a truncated neural network $g^{t_k}\in \mathcal{N}\mathcal{N}_{k}^{t_k}$ satisfying
\begin{equation}
    \|f-g^{t_k}\|_{L^2(\mathbb{B}_1^p,X)}\leq c_1\cdot  k^{-\frac{1}{2}-\alpha'(s,p)}.
\end{equation}
Therefore, we have
\begin{equation}\label{ndjkndfjkasndkj,}
     Conv(\mathbb{P}_{s,1}^p) \subseteq \bigcup_{g^{t_k}\in \mathcal{N}\mathcal{N}^{t_k}_{k}}{\Big\{d\in L^2(\mathbb{B}_1^p,X):\|d-g^{t_k}\|_{L^2(\mathbb{B}_1^p,X)}}\leq c_1\cdot  k^{-\frac{1}{2}-\alpha'(s,d)}\Big\}.
\end{equation}
Next, we need to find a series of balls to cover the above centers lying in $\mathcal{N}\mathcal{N}_{k}^{t_k}$.   Note that by Lemma S.8. in SM the VC dimension of
$$
\mathcal{N}\mathcal{N}_{k}:= \left\{ \sum_{\iota=1}^k{\beta_\iota\cdot \sigma_s(\theta_\iota^T x+b_\iota)}: \beta_\iota\in\mathbb{R}, b_\iota\in [-1,1], \theta\in \Theta^p\right\}
$$
 is at most $O(W\log W)$ with $W=k(p+2)$. Therefore, the VC dimension of truncated neural networks, $\mathcal{N}\mathcal{N}_{k}^{t_k}$, is also bounded by $O(W\log W)$.  According to  Theorem 2.6.7 in \cite{wellner2013weak}, in the space  $L^2(\mathbb{B}_1^p,X)$ we can  at most use $C\cdot \left(\frac{4}{\varepsilon}\right)^{2W \log W}$ balls with radius $\varepsilon\cdot t_k$ to cover  $\mathcal{N}\mathcal{N}_{k}^{t_k}$ where $C>0$ is a universal constant.   Thus, by choosing $\varepsilon=k^{-\frac{1}{2}-\alpha'(s,p)}\cdot \log^{-1} k$, we know
 \begin{equation}\label{zfgdjhasgdkjahd..}
     \mathcal{N}\mathcal{N}_{k}^{t_k}\subseteq \bigcup_{\iota=1}^{N}{  \left\{d\in L^2(\mathbb{B}_1^p,X): \|\xi_\iota-d\|_{L^2(\mathbb{B}_1^p,X)}\leq c_1\cdot  k^{-\frac{1}{2}-\alpha'(s,p)}\right\}},
 \end{equation}
 where $N=k^{k\log k \cdot(\frac{1}{2}+\alpha'(s,p))}$.  Finally, combing \eqref{zfgdjhasgdkjahd..} and \eqref{ndjkndfjkasndkj,} implies that
\begin{equation}\label{odhsakjhdkjhano.}
  Conv(\mathbb{P}_{s,2}^p) \subseteq \bigcup_{\iota=1}^{N}{  \left\{d\in L^2(\mathbb{B}_1^p,X): \|\xi_\iota-d\|_{L^2(\mathbb{B}_1^p,X)}\leq 2c_1\cdot  k^{-\frac{1}{2}-\alpha'(s,p)}\right\}}.
\end{equation}
Meanwhile, we also have
$$
   B_1(\mathbb{P}_{s,2}^p)\subseteq  \bigcup_{\iota=1}^{N}{  \left\{d\in L^2(\mathbb{B}_1^p,X): \|\xi_\iota-d\|_{L^2(\mathbb{B}_1^p,X)}\leq 2c_1\cdot  k^{-\frac{1}{2}-\alpha'(s,p)}\right\}}
$$
because $ B_1(\mathbb{P}_{s,2}^p)= \overline{Conv(\mathbb{P}_{s,2}^p)}$  and the RHS of \eqref{odhsakjhdkjhano.} is closed in $L^2(\mathbb{B}_1^p,X)$.  In conclusion,  the entropy number of $B_1(\mathbb{P}_{s,2}^p)$ satisfies
$$
   \epsilon_{k\log^2 k}(B_1(\mathbb{P}_{s,2}^p))\leq k^{-\frac{1}{2}-\alpha'(s,p)},
$$
which contradicts equation (3.51) in \cite{siegel2022sharp}.  This completes the proof of \eqref{optimal aga rate population bb}.    \hfill\(\Box\)

\vspace{0.3cm}

\textit{Proof of Theorem \ref{PPRpprmodelrate22}.}
The proof is similar to Theorem \ref{EPPRresultppr}, where we boost $\tilde{m}_{k_n}^{aga}$ by using the ensemble technique.   \hfill\(\Box\)

\vspace{0.3cm}



To prove Theorems \ref{EPPRcons}--\ref{EPPRpprmodelrate2b},   we need  a result about convergence of \textit{random AGA in sample version}, which is similar to AGA in sample version but replace $Dic_n^s$ by
 \begin{equation*}
     Dic^{s,q}_n(\A):=\left\{  \sum_{j=1}^{J_n}{a_j\sigma_s(\theta^\top x_{\A}+b_j):  \theta\in\Theta^q, b_j\in[-2,2], \sum_{j=1}^{J_n}|a_j|\leq 1},
     \right\},
\end{equation*}
where $\A$ is the class of selected random indexes with $Card(\A)=q$ during each iteration. Then, we will consider a larger dictionary
$$
 Dic^{s,q}_n=\bigcup_{\A:Card(\A)=q}{Dic^{s,q}_n(\A)},
$$
which is the union of those $ Dic^{s,q}_n(\A)$'s. Unlike the population version of AGA, we will make no assumption on the distribution of $(X_1,\ldots,X_n)$ in its sample version. Therefore, we will see those asymptotic results about AGA in the sample version are  \textit{free from the distribution of $X$}. However, the cost for this generalization is a slightly lower rate compared with that shown in its population version.

\begin{lemma}[{upper bound of random AGA in sample version}]\label{Random Aga sample version corro}
Let $m_k^{aga,q}$ be the output of AGA  after $k$-th iteration where the input vector is $\mathbb{Y}=(Y_1,\ldots,Y_n)^\top$ and $q$ be the number of variables selected during each iteration. Then, we have the error for training data:
\begin{align}
\E_{ \Xi_k}\|\mathbb{Y}-m_k^{aga,q}(X)\|_n^2&\leq \inf_{h\in [\mathfrak{X}_n]_{1}(Dic^{s,q}_{n})}\Big\{\|\mathbb{Y}-h(X)\|_n^2 \nonumber \\
   &+ c(s,p)\left( \|h\|_{[\mathfrak{X}_n]_{1}(Dic^{s,q}_{n})}\right)^2\cdot k^{-1-\frac{2s-1}{2q}}\Big\}\label{sampledhjkz22dbqqq'}
\end{align}
for every iteration $k\in \mathbb{Z}^+$ and some $c(s,p)>0$.
\end{lemma}

\begin{remark}
     If $m(x)=\sum_{\iota=1}^L{m_\iota(x_{\C_\iota})},x\in\B_1^p$, where $\C_\iota\in\{1,\ldots,p\}$  satisfies $Card(\C_\iota)=q$, then approximation rate in \eqref{sampledhjkz22dbqqq'} only depend on $q$ but not $p$.
\end{remark}


\begin{proof}
    Denote by $\delta_k=m-m_k^{aga,q}$ the approximation error in the $k$-th step  for each $k\geq 0$. Let $g_\iota(\theta_\iota^T x_{\A_\iota})=\sum_{j=1}^{J_n}{a_{j,\iota}\sigma_s(\theta_\iota^Tx_{\A_\iota}+b_{j,\iota})}, 1\leq \iota\leq k,$ be the selected functions in the first $k\in\mathbb{Z}^+$ steps respectively.  In the $k$
-th step, we need to project $m$ onto the linear space
$$
   spaceAGA_{k}:=\left\{\sum_{\iota=1}^{k}\sum_{j=1}^{J_n}{a_{j,\iota}\sigma_s(\theta_\iota^Tx_{\A_\iota}+b_{j,\iota})}:   a_{j,\iota}\in\mathbb{R}  \right\},
$$
where $\theta_\iota$ and $b_{j,\iota}$ are fixed. Recall that  $m_k^{aga,q}$ is the corresponding projector.  We have
\begin{equation}\label{HKBHKJBNKJNNqqq}
   \|\delta_k\|^2\leq \|\delta_{k-1}-\langle \delta_{k-1},g_{k}^o\rangle\cdot g_{k}^o\|^2=\|\delta_{k-1}\|^2-|\langle \delta_{k-1},g_{k}^o\rangle|^2
 \end{equation}
for any $g_k^o\in spaceAGA_{k}$ with $\|g_k^o\|=1$.  Now, we remove those linear components of $g_\iota(\theta_\iota^T x_{\A_\iota})\in \mathcal{D}ic_{n}^s, \iota=1, \ldots, k-1$ in $g_k(\theta_k^T x_{\A_k})$ and define
$$
   g_{k}^o:= \frac{g_k(\theta_k^T x_{\A_k})- \mathbb{P}_{spaceAGA_{k-1}}{(g_k(\theta_k^T x_{\A_k}))}}{\|g_k(\theta_k^T x_{\A_k})- \mathbb{P}_{spaceAGA_{k-1}}{(g_k(\theta_k^T x_{\A_k}))}\|},
$$
where $\mathbb{P}_{spaceAGA_{k-1}}(\cdot)$ denotes projector on the linear space ${spaceAGA_{k-1}}$.  Then, \eqref{HKBHKJBNKJNNqqq}  and the fact that $\delta_{k-1}$  is orthogonal to any element in $spaceAGA_{k-1}$ imply that
\begin{equation*}
    \|\delta_k\|^2\leq \|\delta_{k-1}\|^2-\frac{|\langle \delta_{k-1},g_k\rangle|^2}{\|g_k- \mathbb{P}_{spaceAGA_{k-1}}(g_k)\|^2},
\end{equation*}
where the abbreviation $g_{k,\A_k}:= g_k(\theta_k^T X_{\A_k})$.  Note that once $\A_k$ is given, $g_{k,\A_k}:=\sup_{d\in  Dic^{s,q}_n(\A_k)}{|\langle \delta_{k-1}, d\rangle|}$ is not random. Because conditional on $\Xi_{k-1}$, $\delta_{k-1}$ is independent of $\Xi_{k}$. Taking conditional expectation over $\Xi_{k}$ on both sides of \eqref{HKBHKJBNKJNNqqq} leads to
 \begin{equation}\label{Math111}
        \E_{\Xi_k|\Xi_{k-1}}\|\delta_k\|^2\leq \|\delta_{k-1}\|^2- \E_{\Xi_k|\Xi_{k-1}}\left(\frac{|\langle \delta_{k-1},g_{k,\A_k}\rangle|^2}{\|g_{k,\A_k}- \mathbb{P}_{spaceAGA_{k-1}}(g_{k,\A_k})\|^2}\right).
 \end{equation}
Let ${\C}_0= \argmax_{{\C}\subset\{1,\ldots,p\},Card({\C})=q}\sup_{g\in  Dic^{s,q}_n(\C)}|\langle \delta_{k-1},g\rangle|$.  Note that $\C_0$ depends on $\delta_{k-1}$ and is determined once $\Xi_{k-1}$ is given.  However, to simplify the notation,  $\Xi_{k-1}$  is leaved out when we mention $g_{k,\C_0}$.  Because we have one subset $\mathcal{C}\subseteq\{1,\ldots,p\}$  with $Card(\C)=q$ that contains ${\C}_0$, therefore the inequality below follows
\begin{equation}\label{Math131}
\begin{aligned}
     \E_{\Xi_k|\Xi_{k-1}}\left(\frac{|\langle \delta_{k-1},g_{k,\A_k}\rangle|^2}{\|g_{k,\A_k}- \mathbb{P}_{spaceAGA_{k-1}}(g_{k,\A_k})\|^2}\right) \geq \frac{1}{\binom{p}{q}}\frac{\sup_{d\in Dic_n^{s,q}(\C_0)}{|\langle \delta_{k-1},d\rangle|^2}}{\|g_{k,\C_0}- \mathbb{P}_{spaceAGA_{k-1}}(g_{k,\C_0})\|^2}.
\end{aligned}
\end{equation}
The above inequality implies that we need to bound RHS of \eqref{Math131} away from zero.  Let us first analyze its numerator.

Note that for any $h\in Dic_n^{s,q}$,  we  always have
 \begin{align}
     \|\delta_{k-1}\|^2&=\langle \delta_{k-1},h+m-h\rangle=\langle \delta_{k-1},h\rangle+\langle \delta_{k-1},m-h\rangle\nonumber\\
     &\leq \|h\|_{[L^2(\mathbb{B}_1^p,X)]_{1}(Dic^{s,q}_{n})}\cdot\sup_{{\C}\subset\{1,\ldots,p\},Card({\C})=q}\sup_{d\in Dic_n^{s,q}(\C)}|\langle \delta_{k-1},d\rangle|+\|\delta_{k-1}\|\|m-h\|,\label{Math21}
 \end{align}
where in the first equation, we again use the fact that $\delta_{k-1}$ is orthogonal to $m_{k-1}$; and in the last inequality,  we use the fact  that $Dic_n^{s,q}$  lies  in the convex hull of dictionaries  $Dic^{s,q}_n(\C)$, $\C\subseteq \{1,\ldots,p\}$ with $Card(\C)=q$.  By the basic inequality $\|\delta_{k-1}\|\cdot\|m-h\|\leq \frac{1}{2}(\|\delta_{k-1}\|^2+\|m-h\|^2)$, equation \eqref{Math21} implies that
 \begin{equation}\label{Math22}
    \sup_{d\in Dic_n^{s,q}(\C_0)}|\langle \delta_{k-1},d\rangle|\geq \frac{\|\delta_{k-1}\|^2-\|m-h\|^2}{2\cdot\|h\|_{{[L^2(\mathbb{B}_1^p,X)]_{1}(Dic^{s,q}_{n})}}}.
 \end{equation}
This finishes bounding the numerator of \eqref{Math131}.

In conclusion,   \eqref{Math111},  \eqref{Math131} and \eqref{Math22} imply that
  \begin{equation}\label{Math24}
     \E_{\Xi_k|\Xi_{k-1}}\|\delta_k\|^2-\|\delta_{k-1}\|^2\leq -\frac{\max{(\|\delta_{k-1}\|^2-\|m-h\|^2,0})^2}{4\cdot\left(\|h\|_{[L^2(\mathbb{B}_1^p,X)]_{1}(Dic^{s,q}_{n})}\cdot{c(p,q)}\right)^2}\cdot \Omega(g_{k,\C_0}),
 \end{equation}
where $\Omega(g_{k,\C_0}):= \|g_{k,\C_0}- \mathbb{P}_{spaceAGA_{k-1}}(g_{k,\C_0})\|^2$ is a function only w.r.t. $\Xi_{k-1}
$. Taking expectation w.r.t. $\Xi_{k-1}$ on both sides of \eqref{Math24}, by the law of iterated expectation and Jensen's inequality we obtain the following recursive inequality
\begin{equation}\label{Math25}
    \E_{\Xi_{k}}\|\delta_k\|^2\leq  \E_{\Xi_{k-1}}\|\delta_{k-1}\|^2-\frac{ \left(\E_{\Xi_{k-1}}\max{(\|\delta_{k-1}\|^2-\|m-h\|^2,0})\right)^2}{4\cdot\left(\|h\|_{[L^2(\mathbb{B}_1^p,X)]_{1}(Dic^{s,q}_{n})}\cdot{c(p,q)}\right)^2}\cdot \E_{\Xi_{k-1}}\left(\Omega(g_{k,\C_0})\right),
\end{equation}
where $k\in\mathbb{Z}^+$.

 Let $\Delta_\iota=\E_{\Xi_{\iota}} \left( \|\delta_{\iota}\|^2-\|m-h\|^2 \right),\ 1\le \iota\le k$. At this point, our theorem is equivalent to
\begin{equation} \label{equiv}
\Delta_k\leq M\cdot k^{-1-\frac{2s+1}{2q}},
\end{equation}
where $M=4\left(\|h\|_{[L^2(\mathbb{B}_1^p,X)]_{1}(Dic^{s,q}_{n})}\cdot{c(p,q)}\right)^2$.  To prove \eqref{equiv}, consider two scenarios: $\Delta_{k-1}\leq 0$ and $\Delta_{k-1}\ge 0$.  For the first scenario,  due to the projection step of AGA we have $\E_{\Xi_{k}|\Xi_{k-1}}\|\delta_{k}\|^2\leq \|\delta_{k-1}\|^2$ a.s., which implies that
\begin{equation}\label{Theorem2ogalast}
  \E_{\Xi_{k}}(\|\delta_{k}\|^2)=\E_{\Xi_{k-1}}\E_{\Xi_{k}|\Xi_{k-1}}\|\delta_{k}\|^2\leq \E_{\Xi_{k-1}}\|\delta_{k-1}\|^2.
\end{equation}
Therefore,
$$ \Delta_k =  \E_{\Xi_{k}} \left( \|\delta_{k}\|^2-\|m-h\|^2 \right)  \le \E_{\Xi_{k-1}} \left( \|\delta_{k-1}\|^2-\|m-h\|^2 \right) \leq 0$$
and \eqref{equiv} also holds. Consider the second scenario: $\Delta_{k-1}\geq 0$ where our analysis is heavily based on  recursion \eqref{Math25}.  By replacing $\Delta$ with $\Delta/M$ , we can assume $M=1$. Since $\{\Delta_{i}\}_{i=0}^{k-1}$ is a decreasing sequence by a similar argument to \eqref{Theorem2ogalast}, we have $\Delta_{\iota}\geq 0$  for each $0\leq \iota\leq k-1$. Therefore, it follows from \eqref{Math25}  that
\begin{equation}\label{Math26}
      \Delta_\iota\leq  \Delta_{\iota-1}\left(1-\Delta_{\iota-1}\cdot\E_{\Xi_{\iota-1}}\left(\Omega(g_{\iota,\C_0})\right)\right),\ 1 \leq \iota\leq k.
\end{equation}
Taking logarithm on both sides of \eqref{Math26} leads to
$$
\log(\Delta_\iota)\leq \log(\Delta_{\iota-1})- \Delta_{\iota-1} \E_{\Xi_{\iota-1}}\left(\Omega(g_{\iota,\C_0})\right).
$$
Summing above recursions from $1$ to $k$, we have
\begin{equation}\label{BJNKNJKNKLNaaaaaaa}
   \log(\Delta_k)\leq \log(\Delta_{0})- \sum_{\iota=1}^k\Delta_{\iota-1} \E_{\Xi_{\iota-1}}\left(\Omega(g_{\iota,\C_0})\right)\leq \log(\Delta_{0})- \Delta_{k}\cdot\sum_{\iota=1}^k \E_{\Xi_{\iota-1}}\left(\Omega(g_{\iota,\C_0})\right),
\end{equation}
since $\Delta_k\leq \Delta_{\iota-1}$ for all $\iota=1,\ldots,k$.  In fact, we can assume $\Delta_0\leq 1$ and the reason is shown below.  Note that $\E_{\Xi_{0}}\left(\Omega(g_{1,\C_0})\right)$ is bounded away from zero since $Dic^{s,q}_n$ is upper bounded  and $\|g_{\iota,\C_0}- \mathbb{P}_{spaceAGA_{\iota-1}}(g_{\iota,\C_0})\|\leq \|g_k\|$.  Therefore, we can assume $\E_{\Xi_{0}}\left(\Omega(g_{1,\C_0})\right)\geq 1$ in \eqref{Math26} without loss of generality.  Now, if $\Delta_0>1$, \eqref{Math26} implies $\Delta_1\leq 0$, and thus the proof completes.  Therefore,  \eqref{BJNKNJKNKLNaaaaaaa} can be strengthened to
\begin{equation}\label{GYHJGBKJBNKJNK<""}
    \log(\Delta_k)\leq -\Delta_{k}\cdot \E_{\Xi_{k-1}}\left(\sum_{\iota=1}^k\Omega(g_{\iota,\C_0})\right).
\end{equation}
Note that  for any $ 1\leq  \iota\le k$, we have $g_{\iota,\C_0}\in Dic^{s,q}_n:=\bigcup_{\C:Card(\C)=q}{ Dic^{s,q}_n(\C)}$.  By  Lemma S.1 
in \SM\ and Theorem 4 in \cite{siegel2022sharp},  we know the entropy number  $\epsilon_k(Dic^{s,q}_n)\leq c \cdot k^{-\frac{1}{2}-\frac{2s-1}{2q}}$ for some universal $c>0$. Therefore, we have $\sum_{\iota=1}^k\Omega(g_{\iota,\C_0})\ge ck^{1+\frac{2s-1}{q}}$  according to Lemma 1 in \cite{siegel2022sharp}.  Finally,  when $j=k$ from \eqref{GYHJGBKJBNKJNK<""} we have recursive inequalities about $\Delta_j$:
\begin{equation}\label{VYUJHGGBKJBNJKNLKNLaaqc}
     \log(\Delta_j)\leq -c\Delta_{j}\cdot j^{1+\frac{2s+1}{q}}
\end{equation}
provided that  $\Delta_{j-1}>0$.   Since $\Delta_j, j=0,1,\ldots$ is in  decreasing order,  \eqref{VYUJHGGBKJBNJKNLKNLaaqc} also holds for all $j=1,\ldots, k$.  Following arguments on page 7 of \cite{siegel2022optimal},  it is not difficult to finish the proof  from \eqref{VYUJHGGBKJBNJKNLKNLaaqc}.     
\end{proof}

The following two oracle inequalities  will be used in the proof of Theorem \ref{EPPRcons}-- \ref{EPPRpprmodelrate2b}.

\begin{lemma}[{training errors of AGA estimators}]\label{Maintheoremb}
 Suppose
$
\E(e^{c\cdot Y^2})<\infty
$, and
let the sequence $t_n\asymp \ln{n}$.  Then, for any $h_n(x) = \sum_{\iota=1}^{\infty}\sum_{j=1}^{J_n}{a_{n,j,\iota}\sigma_r(\theta_\iota^\top x+b_{n,j,\iota})}$ with $\theta_\iota\in\R^q$  and bounded $m(x),x\in\B_1^p$, we have
\begin{equation} \label{aaa00}
\E_{\mathcal{D}_n, \Xi_{k_n}} (\I(\hat{m}_{k_n,b}(x))) \leq   c\cdot\frac{ (\sum_{\iota=1}^{\infty}\sum_{j=1}^{J_n}{|a_{n,j,\iota}|})^2 }{k_n^{1+\frac{2s-1}{2q}}}
+ \E\left(m(X)-h_n(X)\right)^2+ \omega^1_n.
\end{equation}
Let $A_n:=\{\max_{1\leq i\leq n}|Y_i|\leq t_n\}$ and  $\mathbb{I}(A_n)$ be the indicator of $A_n$. For any  $ \tau \in\mathbb{Z}^+$ (i.e. number of iterations) and $ \nu > \nu_0 $ where $\nu_0>0$ is a constant, we have
\begin{equation}\label{Thesecondmainlemma}
\E_{\mathcal{D}_{n},\Xi_{\infty}}\left( \I(\hat{m}_{k_b^*,b}(x))\cdot \mathbb{I}(A_n)\right)\leq  c\cdot\frac{  (\sum_{\iota=1}^{\infty}\sum_{j=1}^{J_n}{|a_{n,j,\iota}|})^2 }{\tau^{1+\frac{2s-1}{2q}}} + c\cdot \E\left(m(X)-h_n(X)\right)^2+ \omega_n^2,
\end{equation}
where $ \omega^1_n = {c n^{-1} } pk_n J_n\ln(pk_n J_n) \cdot \ln^3 n $ and $ \omega_n^2 = 2\nu n^{-1} \tau (J_n+q)\ln^5{n}$.
\end{lemma}

\begin{proof} Proofs for \eqref{aaa00} and \eqref{Thesecondmainlemma} are similar;  here
we only give the details for the former while that for the latter is given in \SM. Note that  in each run $1\leq b\leq B$,  the output $m_{k_n,b}$ always locates at the function space:
\begin{equation*}
    \mathcal{H}_{k_n,J_n}:= \left\{\sum_{\iota=1}^{k_n}\sum_{j=1}^{J_n}{a_{j,\iota}\sigma_s(\theta_\iota^\top x+b_{j,\iota})} : a_{j,\iota}\in\mathbb{R},\theta_\iota\in\Theta^p,b_{j,\iota}\in [-2,2] \right\}.
\end{equation*}
The proof of Lemma \ref{Maintheoremb} is based on Theorem 1 in \cite{bagirov2009estimation} ({see also Lemma S.2 in \SM}), which shows  that
\begin{equation}\label{Mainlemmaformula2}
    \begin{aligned}
\E_{\mathcal{D}_n, \Xi_{k_n}} (\I(\hat{m}_{k_n,b}(x))) &\leq 2 \E_{\mathcal{D}_{n},\Xi_{k_n}}\left( \| m_{k_n, b}(X)-\mathbb{Y}\|_n^2-\|m(X)-\mathbb{Y}\|_n^2\right)\\
& \ \ \ +\frac{c\ln^2 n}{n}\cdot\sup_{x_1^n}\ln\left(\mathcal{N}(1/(80nt_n),\mathcal{H}^{t_n}_{k_n,J_n},x_1^n)\right),
\end{aligned}
\end{equation}
where $\mathbb{Y}=(Y_1,\ldots,Y_n)^\top$ and $ m_{k_n, b}$ is defined in Section \ref{ssestimator} and  the covering number of $\mathcal{H}^{t_n}_{k_n,J_n}:=\{\max\{\min\{f, t_n\},-t_n\}:f\in \mathcal{H}_{k_n,J_n}\}$ is denoted by $\mathcal{N}(\cdot,\cdot,\cdot)$ . To give an upper bound for the first part on the RHS of \eqref{Mainlemmaformula2}, let
\begin{align*}
  \Delta &=\|m_{k_n,b}(X)-\mathbb{Y}\|_n^2-\|m(X)-\mathbb{Y}\|_n^2 \\
   &= (\|\mathbb{Y}-m_{k_n,b}\|_{n}^2-\|\mathbb{Y}-h_n\|_{n}^2) + (\|\mathbb{Y}-h_n\|_{n}^2-\|\mathbb{Y}-m\|_{n}^2) := \text{I} + \text{II},
\end{align*}
where $\mathbb{Y}=(Y_1,\ldots,Y_n)^\top\in\R^n$. For $\text{I}$,
by Lemma \ref{Random Aga sample version corro} we have
\begin{equation}\label{Sta1}
 \E_{\Xi_{k_n}|\mathcal{D}_{n}}(\text{I})\leq c\cdot\left( \|h_n\|_{[L^2(\mathbb{B}_1^p,X)]_{1}(Dic_n^{r,q})}\right)^2\cdot k_n^{-1-\frac{2s-1}{2q}}.
\end{equation}
Since $\E_{\mathcal{D}_{n}}\left(\|h_n\|^2_{[\mathfrak{X}_n]_1(Dic_n^{r,q})}\right)\leq  (\sum_{\iota=1}^{\infty}\sum_{j=1}^{J_n}{|a_{n,j,\iota}|})^2 $ for each $n$,  it follows by taking expectation w.r.t. $\mathcal{D}_n$ on both sides of \eqref{Sta1} that
\begin{equation}\label{Sta2}
     \E_{\Xi_{k_n}, \mathcal{D}_{n}}(\text{I})\leq c\left( \sum_{\iota=1}^{\infty}\sum_{j=1}^{J_n}{|a_{n,j,\iota}|}\right)^2\cdot k_n^{-1-\frac{2s-1}{2q}}.
\end{equation}
It is easy to see that $\text{II}$ is independent of $\Xi_{k_n}$, thus we further have
\begin{equation}\label{Sta3}
     \E_{\Xi_{k_n}, \mathcal{D}_{n}}(\text{II})=\E\left((Y-h(X))^2- (Y-m(X))^2\right)=\E\left(m(X)-h(X)\right)^2.
\end{equation}
Combination of \eqref{Sta2} and \eqref{Sta3} yields
\begin{equation}\label{Sta4}
     \E_{\Xi_{k_n}, \mathcal{D}_{n}}(\Delta)\leq
     c\left( \sum_{\iota=1}^{\infty}\sum_{j=1}^{J_n}{|a_{n,j,\iota}|}\right)^2\cdot k_n^{-1-\frac{2s-1}{2p}} + \E\left(m(X)-h(X)\right)^2.
\end{equation}
For the second part on the RHS of \eqref{Mainlemmaformula2},  let $VC(\cdot)$ denote the VC dimension of a function class; see \cite{shalev2014understanding} for the definition. With the result about the VC dimension of neural networks given in  \cite{bartlett1998almost},  we have
\begin{equation}\label{myVCaa}
 VC({\mathcal{H}^{t_n}_{k_n,J_n}}) \leq VC(\mathcal{H}_{k_n,J_n}) \leq c\cdot k_nJ_n(p+1)\ln(k_nJ_n(p+1))
\end{equation}
for some $c>0$. Then, by Lemma 9.2 in \cite{gyorfi2002distribution} we have
\begin{align}\label{Sta5}
    \mathcal{N}(1/(80n t_n),\mathcal{H}^{t_n}_{k_n,J_n},x_1^n) &\leq 3\left( \frac{4e t_n}{1/(80n t_n)}\ln\left(\frac{6e t_n}{1/(80n t_n)} \right) \right)^{VC(\mathcal{H}^{t_n}_{k_n,J_n})} \\
    & \leq 3\left( \frac{4e t_n}{1/(80n t_n)}\ln\left(\frac{6e t_n}{1/(80n t_n)} \right) \right)^{VC(\mathcal{H}_{k_n,J_n})}\nonumber \\
    & \leq 3\left( \frac{4e t_n}{1/(80n t_n)}\ln\left(\frac{6e t_n}{1/(80n t_n)} \right) \right)^{c\cdot pk_nJ_n\ln(pk_nJ_n)}\nonumber \\
    & \leq 3\left( 480e n t_n^2 \right)^{c\cdot k_nJ_n \ln(k_nJ_n)}.\nonumber
\end{align}
It follows from \eqref{Sta5} and \eqref{Mainlemmaformula2} that
    \begin{align}\label{Sta6}
      \E_{\mathcal{D}_n, \Xi_{k_n}} (\I(\hat{m}_{k_n,b}(x)))  &\leq  \frac{c\ln^2 n}{n}pk_nJ_n\ln(pk_nJ_n) \ln\left(n t_n \right) +\E_{\Xi_{k_n}, \mathcal{D}_{n}}(\Delta).
    \end{align}
Therefore, Lemma \ref{Maintheoremb} follows from \eqref{Sta6} and   \eqref{Sta4}.
\end{proof}


\textit{Proof of Theorem  \ref{EPPRcons}.}\label{Theorem 666}
Without loss of generality, we can assume $L=1$. Since the class of infinitely differentiable functions $C^\infty(\B_1^q)$  is dense in $L^2(\B_1^q,X)$, given any $\epsilon>0$ one can find $\tilde{m}(x)\in C^\infty(\B_1^q)$ such that
\begin{equation}\label{JKHNKJHNwera}
    \E\left( m(X)-\tilde{m}(X)\right)^2 \leq \epsilon/4.
\end{equation}
According to \cite{LESHNO1993861}, there is  a series of functions $h_n(x) = \sum_{\iota=1}^{k_n}\sum_{j=1}^{J_n}{a_{n,j,\iota}\sigma_r(\theta_\iota^\top x+b_{n,j,\iota})}$ with $\theta_\iota\in\R^q$ and $n_1(\epsilon)\in\mathbb{Z}^+$,  such that
\begin{equation}
   \E\left( \tilde{m}(X)-h_n(X)\right)^2 \leq  \epsilon/4
\end{equation}
when $n\ge n_1(\epsilon)$ and $J_n\to\infty$.  As a consequence,   there is  $h_{n_1(\epsilon)}\in \mathcal{H}_{k_n,J_n}$ such that
\begin{equation*}\label{dfngmdfm}
  \text{I}:= \E(m(X) - h_{n_1(\epsilon)}(X))^2 \le \epsilon/2
\end{equation*}
no matter the choice of $X$. Since  $\mathcal{H}_{k_n,J_n}$ is in increasing order,  we know
$h_{n_1(\epsilon)}\in \mathcal{H}_{k_n,J_n}$ for each $n\ge n_1(\epsilon)$. Thus, it is valid to replace $h_n$ in \eqref{aaa00} with $h_{n_1(\epsilon)}$ when $n\ge n_1(\epsilon)$. On the other hand, as $ k_n \to \infty $ there exits $n_2(\epsilon)\in\mathbb{Z}^+$ such that for any $n\ge n_2(\epsilon)$, we  have
$$
\text{II}:= c\cdot \left( \sum_{\iota=1}^{k_{n_1(\epsilon)}}\sum_{j=1}^{J_{n_1(\epsilon)}}{|a_{n,j,\iota}|}\right)^2 \cdot {k_n^{-1-\frac{2s-1}{2p}}} \le \epsilon/2.
$$
By \eqref{aaa00} of Lemma \ref{Maintheoremb}, there exists $n(\epsilon) \ge \max(n_1(\epsilon), n_2(\epsilon))$ such that for any $n>n(\epsilon)$,
\begin{equation*}
\MSE(\hat{m}_{k_n,b}(x)) \leq   \text{I} + \text{II} + \omega^1_n \le \epsilon,
\end{equation*}
which implies
\begin{equation} \label{dfknge}
\MSE(\hat{m}_{k_n,b}(x)) \to 0,\ \text{as}\ n\to\infty.
\end{equation}
Note that, by Jensen's inequality,
\begin{equation}\label{JensenAB}
    |\hat{m}^{ens}_{k_n}(x)-m(x)|^2\leq \frac{1}{B}\sum_{b=1}^B{| \hat{m}_{k_n,b}(x) - m(x)|^2}.
\end{equation}
Therefore,
\begin{equation}\label{JensenAB1}
    \MSE(\hat{m}^{ens}_{k_n}(x))\leq \frac{1}{B}\sum_{b=1}^B \MSE(\hat{m}_{k_n,b}(x)).
\end{equation}
As a consequence, Theorem  \ref{EPPRcons} follows immediately from both \eqref{dfknge} and \eqref{JensenAB1}.  \hfill\(\Box\)

\vspace{0.3cm}

\textit{Proof of Theorem  \ref{EPPRresultppr}.}
Recall that
$
  m(x)=\sum_{\iota=1}^\infty {m_\iota(\theta_\iota^\top  x_{\C_\iota})},
$
where $\theta_\iota\in\mathbb{R}^{q}$ and $ card(\C_\iota)=  q $. By  Lemma 
S.6 in \SM, it is not difficult to derive a stronger result as follows. In fact, there is a sequence of functions  $$h_n(x)= \sum_{\iota=1}^\infty\sum_{j=1}^{J_n}{a_{n,j,\iota}\sigma_s(\theta_\iota^\top x+b_{n,j,\iota})} \in c(s,m)\cdot B_1(Dic_n^{r,q})$$  with $\sum_{\iota=1}^{\infty}\sum_{j=1}^{J_n}{|a_{n,j,\iota}|}<\infty$  and  $|b_{n,j,\iota}|\leq 2$,  such that
\begin{equation*}
        \sup_{x\in [-1,1]}{|m(x)- h_n(x)|}\leq c(s)\sum_{\iota=1}^\infty{TV(m_\iota^{(s)})}J_n^{-s},
\end{equation*}
where $c(s,m)>0$ does not depend on $n\in\mathbb{Z}^+$ while $c(s)$ only relates to $s$.

By \eqref{aaa00} of Lemma \ref{Maintheoremb},  we immediately have
\begin{equation}\label{Theorem4part1}
\MSE(\hat{m}_{k_n,b}(x)) \leq   c \cdot k_n^{-1-\frac{2s-1}{2q}}
+ c\cdot J_n^{-2s} + \omega^1_n.
\end{equation}
Properly selecting $k_n $ and $ J_n$ such that $k_n^{-1-\frac{2s-1}{2p}}=J_n^{-2s}=\frac{k_nJ_n}{n}$, namely $k_n=n^{(1+\frac{2s-1}{2q})/(2s)}$, the RHS of  \eqref{Theorem4part1} is bounded by $c\cdot n^{-\alpha_*(s,q)}\cdot\ln^4 n$, where
$$\alpha_*(s,q)= (\frac{s+.5}{q}+1)/\left((\frac{s+.5}{q}+1)+ (1+\frac{1}{2s})(1-\frac{1}{2q})\right).$$
Finally,  Theorem  \ref{EPPRresultppr} follows from  \eqref{Theorem4part1} and \eqref{JensenAB}.  \hfill\(\Box\)

\def\I{\mathcal{I}}

\vspace{0.3cm}

\textit{Proof of Theorem \ref{EPPRcons2a}.}
Without loss of generality, we assume $L=1$. Recalling $A_n $ in Lemma \ref{Maintheoremb}, let $ A_n^c $ be the complement. Write
\begin{align}
    \I(\hat{m}_{k_b^*,b}(x))
    &=\I(\hat{m}_{k_b^*,b}(x))\cdot \mathbb{I}(A_n^c)+ \I(\hat{m}_{k_b^*,b}(x)))\cdot \mathbb{I}(A_n) := \text{I}+\text{II}. \nonumber
\end{align}
For any given $\epsilon>0$, we have
\begin{equation} \label{Theorem5part2}
\begin{aligned}
\P\left( \I(\hat{m}_{k_b^*,b}(x))>\epsilon\right)&\leq \P(\text{I}>\epsilon/2)+\P(\text{II}>\epsilon/2)\\
   &\leq \P(A_n^c)+\P(\text{II}>\epsilon/2):=\text{III}+\text{IV}.
\end{aligned}
\end{equation}

First, we consider part III.  By the distribution assumption on $ Y $, we have
\begin{align}
\label{subgaussguanjian}    \P(A_n^c)&=1- \P\left(\max_{1\leq i\leq n}|Y_i|\leq t_n\right)
    =1-\left[\P(|Y_i|\leq t_n)\right]^n
    \leq 1-(1-c\cdot e^{-c\cdot t_n^2})^n \\
    &\le 1-e^{n\cdot \ln(1-c\cdot e^{-c\cdot t_n^2})}. \nonumber
\end{align}
Therefore,
$$
 \text{III} = \P(A_n^c) \to 0.
$$

Next, consider part IV using \eqref{Thesecondmainlemma} in Lemma \ref{Maintheoremb}.
Choose the same  $h_n$ in the proof of Theorem \ref{EPPRcons}. For any given $\eta>0$, there is $n_1(\eta)\in\mathbb{Z}^+$ such that
\begin{equation}\label{adfgarohr}
  \text{V}:=\E\left(m(X)-h_{n_1(\eta)}(X)\right)^2 \leq c\cdot J_{n_1(\eta)}^{-2s} \le \eta/3.
\end{equation}
Now, fix $ n_1(\eta) $ for the given $ \eta $ and replace $h_n$ in \eqref{Thesecondmainlemma} with $h_{n_1(\eta)}$ for each $n\geq n_1(\eta)$. Choose $\tau_n=\lfloor \ln{n} \rfloor$ which goes to infinity as $n\to\infty$.
Therefore, there is $n_2(\eta)\in\mathbb{Z}^+$ such that for any $n\ge n_2(\eta)$
\begin{equation}\label{adfgarohr2}
    \text{VI}:=  c\cdot \tau_n^{-1-\frac{2s-1}{2q}}\cdot\left(\sum_{\iota=1}^{k_{n_1(\eta)}}\sum_{j=1}^{J_{n_1(\eta)}}{|a_{n,j,\iota}|}\right)^2  + \omega_n^2\le \eta/3.
\end{equation}
Combination of \eqref{adfgarohr} and \eqref{adfgarohr2} leads to
\begin{equation*}
    \E_{\mathcal{D}_{n},\Xi_{\infty}}\left( \I(\hat{m}_{k_b^*,b}(x)) \cdot \mathbb{I}(A_n)\right) \leq  \text{V}+\text{VI}\leq \eta
\end{equation*}
for any $n\ge \max\{n_1(\eta),n_2(\eta)\}$, which implies that for every $ b $,
 \begin{equation}\label{Theorem5lastone}
     \lim_{n\to\infty}{\E_{\mathcal{D}_{n},\Xi_{\infty}} \left(\I(\hat{m}_{k_b^*,b}(x))\cdot\mathcal{I}(A_n)\right)}=0.
 \end{equation}
By  Chebyshev inequality,  above equation \eqref{Theorem5lastone} implies that
$$
\text{IV}=\P(\text{II}>\epsilon/2)\to 0.
$$

In conclusion, we have proved that the LHS of \eqref{Theorem5part2} converges to $0$. On the other hand, by Jensen's inequality, we have
\begin{equation}\label{JensenCD}
    |\hat{m}^{ens}_{*}(x)-m(x)|^2\leq \frac{1}{B}\sum_{b=1}^B{| \hat{m}_{k_b^*,b}(x) - m(x)|^2}.
\end{equation}
Finally, Theorem \ref{EPPRcons2a} follows from  \eqref{Theorem5lastone} and \eqref{JensenCD}.  This completes the proof.  \hfill\(\Box\)

\vspace{0.3cm}

\textit{Proof of Theorem \ref{EPPRpprmodelrate2b}.}
For each $n\in\mathbb{Z}^+$, let $A_n$ and $A_n^c $ be defined in the proof of Theorem \ref{EPPRcons2a} and $a_n=n^{-\alpha_*(s,q)}\cdot \ln^5 n$, where $\alpha_*(r,q)$ is given in \textit{Proof of Theorem \ref{EPPRresultppr}}. Write
\begin{equation*}
    \frac{1}{a_n}\cdot\I(\hat{m}_{k_b^*,b}(x))=\frac{1}{a_n}\cdot\I(\hat{m}_{k_b^*,b}(x)\cdot \mathbb{I}(A_n^c)
    + \frac{1}{a_n}\cdot\I(\hat{m}_{k_b^*,b}(x))\cdot \mathbb{I}(A_n) := \text{I}+\text{II}.
\end{equation*}
For any given $\delta>0$  and $M(\delta)>0$, we have
\begin{equation}\label{Theorem7part2}
\begin{aligned}
   \P\left(\frac{1}{a_n}\cdot\I(\hat{m}_{k_b^*,b}(x))>M(\delta)\right)&\leq\P(\text{I}>M(\delta)/2)
   +\P(\text{II}>M(\delta)/2)\\
   &\leq \P(A_n^c)+\P(\text{II}>M(\delta)/2):=\text{III}+\text{IV}.
\end{aligned}
\end{equation}
For III, we can use the same argument used in \textit{Proof of Theorem  \ref{EPPRcons2a}} to show that there {is $n_1(\delta)\in\mathbb{Z}^+$} such that  $\text{III} < \delta/2$ for all $n\geq n_1(\delta)$.

Next, we consider part IV. Choose the same $h_n$ in \textit{Proof of Theorem  \ref{EPPRresultppr}} satisfying
\begin{equation*}
        \sup_{x\in [-1,1]}{|m(x)- h_n(x)|}\leq c(s)\sum_{\iota=1}^\infty{TV(m_\iota^{(s)})}J_n^{-s}
\end{equation*}
for  some $c(s)>0$.  By \eqref{Thesecondmainlemma} of Lemma \ref{Maintheoremb} and similar arguments in \textit{Proof of Theorem \ref{EPPRresultppr}}, we have
\begin{equation}\label{Theorem7part1}
\E_{\mathcal{D}_{n},\Xi_{\infty}}\left( \I(\hat{m}_{k_b^*,b}(x)\cdot \mathbb{I}(A_n)\right)\leq   c\cdot \tau_n^{-1-\frac{2s-1}{2q}}
+ c\cdot J_n^{-2s}+ \omega_n^2
\end{equation}
for any $\tau_n \in\mathbb{Z}^+$.  Properly selecting $\tau_n$ and $J_n$ such that $\tau_n^{-1-\frac{2s-1}{2q}}=J_n^{-2s}= n^{-1} \tau_n J_n$, namely $\tau_n=n^{(1+\frac{2s-1}{2q})/(2s)}\to\infty$, then the LHS of  \eqref{Theorem7part1} is bounded by the product of a positive constant and $a_n=n^{-\alpha_*(s,q)}\cdot \ln^5 n$. The application of Chebyshev inequality implies that there exists a large $M_1(\delta)>0$ such that
$$
   \text{IV}:=\P(\text{II}>M_1(\delta)/2)\leq \frac{2\cdot \E_{\mathcal{D}_{n},\Xi_{\infty}}(\text{II})}{M_1(\delta)}\leq \delta/2
$$
holds for all $n\in\mathbb{Z}^+$.

Therefore, by \eqref{Theorem7part2} we have shown that
\begin{equation}\label{Aldfng}
 \P\left(\frac{1}{a_n}\cdot\I(\hat{m}_{k_b^*,b}(x))>M_1(\delta)\right)\le \text{III}+\text{IV}\le \delta
\end{equation}
holds for each $n\ge n_1(\delta)$. Theorem \ref{EPPRpprmodelrate2b} follows from
  \eqref{Aldfng} and \eqref{JensenCD}.  \hfill\(\Box\)

\section{Discussion} \label{secConculsion}

As commented by \cite{hastie01}, ``PPR does represent an important intellectual advance, one that has blossomed in its reincarnation in the field of neural networks''. \cite{hastie01} attributed the low popularity of PPR to the computational issue, which is part of the motivation for this paper. This paper studied three greedy algorithms of PPR and proposed an ensemble procedure, ePPR, based on ``feature boosting''. Our extensive numerical study suggests that ePPR can generally produce more accurate predictions than the popular methods, including RF and its variants, for data with either continuous response or categorical responses.

One possible reason for the better performance of ePPR over RF is their underlying models, PPR and CART. The former can approximate arbitrary functions, as shown in this paper, while the latter cannot because its splitting is based on marginal variables.  Moreover, PPR is continuous but CART is a step function which may cause problems in prediction (\cite{wager2015adaptive}). The main problem with ePPR is the computational burden. However, as an ensemble procedure, the calculation of ePPR can be easily reduced by using parallel computation.

Our work also sheds light on the long debate about when ANN and RF have advantages over the other. As a special case of ANN, the performance of ePPR suggests that when statistical tuning is incorporated, ANN can perform as well as or even better than RF for data of small to medium size, challenging the commonly received understanding  that ANN can only have  advantages when the sample size is large.


\bigskip

\begin{center}
SUPPLEMENTARY MATERIAL
\end{center}
\textbf{
Supplement to: "Ensemble Projection Pursuit  for General Nonparametric Regression"}
provides some auxiliaries used in the proofs in Section \ref{SecProof} and other relevant results.

\bigskip

\baselineskip1.5em
\bibliographystyle{imsart-number} 
	\bibliography{ref}

\end{document}